\pdfoutput=1
\documentclass[a4paper,12pt]{article}
\usepackage{amssymb}
\usepackage{amsmath}
\usepackage{amscd}
\usepackage[arrow,matrix]{xy}
\usepackage[bookmarksnumbered,colorlinks,linkcolor=black,citecolor=black,urlcolor=black]{hyperref}

\newtheorem{theo}{Theorem}[section]
\newtheorem{prop}[theo]{Proposition}
\newtheorem{lem}[theo]{Lemma}
\newtheorem{cor}[theo]{Corollary}
\newtheorem{defi}[theo]{Definition}
\newtheorem{rem}[theo]{Remark}

\def \K{{\mathcal K}}

\def \Br {{\rm{Br}}}

\def \Ga {{\Gamma}}
\def \R {{\mathbb{R}}}
\def \Pic {{\rm {Pic}}}

\def \Gal {{\rm{Gal}}}

\def \Im {{\rm {Im}}}

\def \A{{\mathbb A}}
\def \P{{\mathbb P}}
\def \S{{\mathbb S}}
\def \K{{\mathbb K}}

\def \dim {{\rm{dim}}}
\def \Hom {{\rm {Hom}}}
\def \End {{\rm {End}}}
\def \Pic {{\rm {Pic}}}
\def \GL {{\rm {GL}}}

\def \SO {{\rm {SO}}}
\def \bSO {{\bf {SO}}}

\def \Aut{{\rm Aut}}

\def \Sh {{\rm {Sh}}}
\def \tr {{\rm {tr}}}

\def\ov{\overline}

\def \Z {{\mathbb Z}}
\def \Q {{\mathbb Q}}
\def \F {{\mathbb F}}

\def \Id {{\rm Id}}

\def \Tr {{\rm{Tr}}}
\def \Mat {{\rm{Mat}}}
\def \sym {{\rm{sym}}}

\def \rk {{\rm{rk}}}

\def\G{{\mathbb G}}

\def\C{{\mathbb C}}

\def\lra{\longrightarrow}
\def\cl{{\rm cl}}
\def\N{{\rm N}}
\def\H{{\rm H}}

\def\Tr{{\rm Tr}}

\def\NS{{\rm NS\,}}

\def\cl{{\rm cl}}

\def\discr{{\rm discr}}
\def\Ga{\Gamma}

\def\et{{\rm{\acute et}}}

\def\exp{{\rm exp}}

\newcommand{\Xcld}{X^{\mathrm{cl},\leq d}}
\newcommand{\MT}{\mathrm{MT}}

\newcommand{\bthe}{\begin{theo}}
\newcommand{\ble}{\begin{lem}}
\newcommand{\bpr}{\begin{prop}}
\newcommand{\bco}{\begin{cor}}
\newcommand{\bde}{\begin{defi}}
\newcommand{\ethe}{\end{theo}}
\newcommand{\ele}{\end{lem}}
\newcommand{\epr}{\end{prop}}
\newcommand{\eco}{\end{cor}}
\newcommand{\ede}{\end{defi}}
\newcommand{\brem}{\begin{rem}}
\newcommand{\erem}{\end{rem}}

\usepackage[textheight=220mm,textwidth=150mm]{geometry}

\title{On uniformity conjectures for abelian varieties\\ and K3 surfaces}
\author{Martin Orr, Alexei N. Skorobogatov and Yuri G. Zarhin}
\date{\today}
\begin{document}
\maketitle

\begin{abstract}
\noindent We discuss logical links among uniformity conjectures concerning K3 surfaces and abelian
varieties of bounded dimension defined over
number fields of bounded degree.
The conjectures concern the endomorphism algebra of an abelian variety,
the N\'eron--Severi lattice of a K3 surface, and the Galois invariant subgroup of the geometric
Brauer group.
\end{abstract}

{\small
\tableofcontents
}

\footnotetext{The first and second named authors have been supported by the EPSRC grant EP/M020266/1.
The third named author is partially supported by Simons Foundation Collaboration grant \#~585711.
The second and third named authors would like to thank the Max Planck Institut f\"ur Mathematik in Bonn
for hospitality and support. We are grateful to the organisers of the workshop 
``Arithmetic of curves" at Baskerville Hall for excellent working conditions that enabled us
to complete this project.}

\parskip=2pt minus1pt
\parindent=12pt
\baselineskip=15pt
\abovedisplayskip=8pt plus4pt minus4pt
\belowdisplayskip=8pt plus4pt minus4pt

\section{Introduction}

The aim of this paper is to explore logical links among several
conjectures about K3 surfaces and abelian varieties defined over number fields.
These conjectures state that certain
invariants take only finitely many values provided the degree of the field of definition
and the dimension (in the case of abelian varieties) are bounded.

Let $k$ be a number field with algebraic closure $\bar k$ and let $\Ga=\Gal(\bar k/k)$.
For a variety $X$ over $k$ we write $\ov X=X\times_k\bar k$.

\medskip

\noindent{\bf Coleman's conjecture about $\End(\ov A)$.}
{\em Let $d$ and $g$ be positive integers.
Consider all abelian varieties $A$ of dimension $g$ defined over number fields of degree $d$.
Then there are only finitely many isomorphism classes among the rings $\End(\ov A)$.}

\medskip

This or a closely related conjecture is attributed to Robert Coleman in 
\cite[Remark 4]{Sha96}; see also
Conjecture C$(e,g)$ in \cite[p.~384]{BFGR} and the second paragraph of \cite[p.~651]{MW94}. There is
a version of this conjecture in which $\End(\ov A)$
is replaced by the ring $\End(A)$ of endomorphisms of $A$ defined over $k$.
It is not too hard to show that Coleman's conjecture
about $\End(\ov A)$ is equivalent to Coleman's conjecture about $\End(A)$, see Theorem \ref{22a1}.

In his recent paper R\'emond proved that Coleman's conjecture implies the uniform boundedness
of torsion $A(k)_{\rm tors}$ and of the minimal degree of an isogeny between isogenous abelian
varieties, see \cite[Thm.~1.1]{Rem18}. In this paper we would like to point out several other
consequences of Coleman's conjecture.

\medskip

\noindent{\bf Shafarevich's conjecture about $\NS(\ov X)$.} 
{\em Let $d$ be a positive integer.
There are only finitely many lattices $L$, 
up to isomorphism, for which there exists a K3 surface $X$ defined over
a number field of degree $d$ such that $\NS(\ov X)\cong L$.}

\medskip

It is in this form that Shafarevich has stated his conjecture in \cite{Sha96}. 
Since there are only finitely
many lattices of bounded rank and discriminant \cite[Ch.~9, Thm.~1.1]{cassels:quadratic-forms}, Shafarevich's conjecture is equivalent
to the boundedness of the discriminant of $\NS(\ov X)$.
One can also state a variant of Shafarevich's conjecture in which $\NS(\ov X)$
is replaced by its Galois-invariant subgroup $\NS(\ov X)^\Ga$, or, alternatively, by $\Pic(X)$.
In Theorem \ref{Sh-eq}
we show that all these versions of Shafarevich's conjecture are equivalent.

Similarly to Shafarevich's conjecture, Coleman's conjecture can be restated in terms of lattices. 
Recall that $\End(A)$ is an order in the semisimple $\Q$-algebra $\End(A)_\Q=\End(A)\otimes\Q$.
Let us define $\discr(A)$ as the discriminant of the integral
symmetric bilinear form $\tr(xy)$ on $\End(A)$,
where $\tr:\End(A)_\Q\to\Q$ is the reduced trace.
An equivalent form of Coleman's conjecture says
that $\discr(A)$ is uniformly bounded for abelian varieties $A$
of bounded dimension defined over number fields of bounded degree.

\medskip

We denote by $\Br(X)=\H^2_\et(X,\G_m)$ the (cohomological) Brauer group of a 
scheme $X$. When $X$ is a variety over a field $k$, we use the standard notation $\Br_0(X)$
for the image of the canonical map $\Br(k)\to \Br(X)$ and $\Br_1(X)$ for the kernel of the map $\Br(X) \to \Br(\ov X)^\Ga$.

Assume that $k$ is finitely generated over $\Q$,
for example, $k$ is a number field.
The geometric Brauer group $\Br(\ov X)$ has a natural structure of a $\Ga$-module. By the main result of \cite{SZ08}, if $X$ is an abelian variety or a K3 surface over~$k$, then $\Br(\ov X)^\Ga$ is finite.
Furthermore, $\Br_1(X) / \Br_0(X)$ injects into $\H^1(k, \Pic(\ov X))$ and if $X$ is a K3 surface over $k$, then $\H^1(k, \Pic(\ov X))$ is finite \cite[Remark~1.3]{SZ08}.
Consequently the finiteness of $\Br(X)/\Br_1(X) \subset \Br(\ov X)^\Ga$ implies that $\Br(X)/\Br_0(X)$ is finite when $X$ is a K3 surface over a field finitely generated over $\Q$.
It is well known that the cardinality of the finite group $\H^1(k, \Pic(\ov X))$, where $X$ is a K3 surface over an arbitrary field $k$ of characteristic zero, is bounded, see e.g.\ \cite[Lemma 6.4]{VAV17}.

\medskip

\noindent{\bf V\'arilly-Alvarado's conjecture.} \cite[Conj. 4.6]{VA17}
{\em Let $d$ be a positive integer and let $L$ be a primitive sublattice of the K3 lattice $E_8(-1)^{\oplus 2} \oplus U^{\oplus 3}$.
If $X$ is a K3 surface defined over
a number field of degree $d$ such that $\NS(\ov X)\cong L$, 
then the cardinality of $\Br(X)/\Br_0(X)$ is bounded.}

\medskip

A stronger form of this conjecture omits the reference to the N\'eron--Severi lattice. It concerns 
the uniform boundedness of the Galois invariant subgroup of the geometric Brauer group.

\medskip

\noindent {\bf Conjecture Br(K3).} {\em Let $d$ be a positive integer.
There is a constant $C = C(d)$ such that, if $X$ is a K3 surface defined over a number field of degree $d$, then $|\Br(\ov X)^\Ga| < C$.}

\medskip

A similar conjecture can be stated for abelian varieties of given dimension.

\medskip

\noindent {\bf Conjecture Br(AV).} 
{\em Let $d$ and $g$ be positive integers.
There is a constant $C = C(d,g)$ such that, if $A$ is an abelian variety of dimension~$g$ defined over
a number field of degree~$d$, then $|\Br(\ov A)^\Ga| < C$.}

\medskip

The main results of this paper are summarised in the following diagram:
$$\left. \begin{array}{ccc}
\text{Coleman's conjecture}&\Longrightarrow&\text{Shafarevich's conjecture}\\
\Downarrow&&\\
\text{Br(AV)}&\Longrightarrow &\text{V\'arilly-Alvarado's conjecture}
\end{array}\right\}\Longrightarrow \text{Br(K3)}$$

We shall now discuss some known results in the direction of these conjectures.
All of the aforementioned conjectures hold for abelian varieties and K3 surfaces
with complex multiplication \cite{OS}.
Coleman's conjecture for elliptic curves follows from the Brauer--Siegel theorem.
Fit\'e and Guitart \cite{FG} have recently made progress on Coleman's conjecture in the case $g=2$, $d=1$ by showing that there are only finitely many $\Q$-algebras which can appear as $\End(\ov A) \otimes \Q$ when $A/\Q$ is an abelian surface such that $\ov A$ is not simple.
Shafarevich \cite{Sha96b} proved a function field version of Coleman's conjecture for abelian surfaces whose endomorphism algebra is a quaternion algebra.

All these conjectures may be stated in the form ``in a certain class of moduli spaces, only finitely many spaces in the class have rational points over number fields of degree~$d$, excluding points which lie in subvarieties of positive codimension parameterising objects with extra structures.''
From this point of view, Nadel~\cite{Nad89} and Noguchi \cite{Nog91} have proved similar conjectures over complex function fields, while Abramovich and V\'arilly-Alvarado \cite{AVA18} have proved that Lang's conjecture implies a result of this type for abelian varieties with full level structure.

Some of the above conjectures are known for the fibres of one-parameter families.
In particular, a result of Cadoret and Tamagawa \cite{CT13} implies Coleman's conjecture within a one-parameter family of abelian varieties (see Appendix).
Cadoret and Charles \cite{CC} have proved uniform boundedness of the \( \ell \)-primary subgroup of the Brauer group for one-parameter families of abelian varieties and K3 surfaces.
V\'arilly-Alvarado and Viray obtained bounds for the Brauer group for one-parameter families of Kummer surfaces attached to products of isogenous elliptic curves [VAV17, Thm.~1.8].

\medskip

Here is an outline of the paper.
After discussing some preliminary results in Section \ref{1}, we establish the equivalence of various
forms of Coleman's conjecture and also those of Shafarevich's conjecture in Section \ref{2}.

Section \ref{3} is devoted to proving that Coleman's conjecture implies Br(AV). We give two
different proofs that uniformly large primes do not divide $|\Br(\ov A)^\Ga|$.
In Section \ref{Sect. 3.1} we give a shorter proof based on the aforementioned result of R\'emond \cite[Thm.~1.1]{Rem18} and the methods of \cite{Zar77, Zar85}. 
In Section \ref{3.b}
we give a proof that does not use \cite[Thm.~1.1]{Rem18}; this approach has the advantage 
of being more general as it applies also to finitely generated fields. 
Here the key role is played by 
the image $\Lambda_\ell(A)$ of the $\ell$-adic group algebra of the Galois group in the endomorphism ring
of the $\ell$-adic Tate module $T_\ell(A)$. A crucial observation (Theorem \ref{autre1}) is that
a matrix algebra over the opposite algebra of
$\Lambda_\ell(A)$ is isomorphic to $\End(B)\otimes\Z_\ell$, where $B$ 
is an abelian variety isogenous to an abelian subvariety of a bounded power of $A$. 
Hence $\discr(\Lambda_\ell(A))$ divides $\discr(B)$, so
under Coleman's conjecture we obtain an upper bound for $\discr(\Lambda_\ell(A))$.
The relevance of this to Br(AV) is that a prime $\ell>4\dim(A)$ dividing $|\Br(\ov A)^\Ga|$
must also divide $\discr(\Lambda_\ell(A\times A^\vee))$, where $A^\vee$
is the dual abelian variety of $A$. To complete the proof that Coleman's conjecture implies Br(AV)
one needs to show the uniform boundedness of the $\ell$-primary torsion of $\Br(\ov A)^\Ga$ 
for a fixed $\ell$; this proof can be found in Section \ref{Sect. 3.2}.
In Section~\ref{converse}, we prove some partial converses to Theorem~\ref{3.2}: bounds for Brauer groups of abelian varieties imply information about their endomorphisms.

In Section \ref{5} we use the K3 surfaces version of Zarhin's trick
from \cite{OS} to produce a uniform Kuga--Satake construction that does not
depend on the degree of polarisation. The Hodge-theoretic aspect of this construction allows us
to show that Coleman's conjecture implies Shafarevich's conjecture.
In Section \ref{6} we use the compatibility with Galois action
to prove that Br(AV) implies V\'arilly--Alvarado's conjecture.
By the finiteness of the isomorphism classes of lattices of the same rank and discriminant,
it is clear that the conjectures of Shafarevich and V\'arilly--Alvarado together imply 
Conjecture Br(K3).

In the appendix to this paper we deduce the conjectures of Coleman and Shafarevich for one-parameter families from results of Cadoret--Tamagawa and Hui.

\section{Preliminaries} \label{1}

\subsection{Lattices}

In this paper we refer to a free abelian group $L$ of finite positive rank with
a non-degenerate integral symmetric bilinear form $(x.y)$ as a {\em lattice}.
Write $L^*=\Hom(L,\Z)$ and $L_\Q=L\otimes_\Z\Q$.
The {\em discriminant group} of a lattice $L$ is defined as the cokernel of the map
$L\to L^*$ sending $x\in L$ to the linear form $(x.y)$. 
The {\em discriminant} $\discr(L)$ of $L$ is the determinant
of the matrix $(e_i.e_j)$, where $e_1,\ldots, e_n$ is a $\Z$-basis of $L$.
This is independent of the choice of basis $e_1, \ldots, e_n$.
We have $|\discr(L)|=|L^*/L|$.

Let $\ell$ be a prime.
We define the {\em discriminant} $\discr(L)$ of a free $\Z_\ell$-module $L$ of finite positive rank
equipped with a symmetric $\Z_\ell$-valued bilinear form in the same way.
However, in this case there is an ambiguity coming from the choice of $\Z_\ell$-basis for~$L$:
$\discr(L)$ is well-defined up to multiplication by a square in $\Z_\ell^\times$.
In practice, every use we make of the discriminant of a $\Z_\ell$-module $L$ will only depend on the $\ell$-adic valuation of $\discr(L)$, which is well-defined.

\ble \label{2.1}
Let $L$ be a lattice with discriminant $d$. Let $G$ be a finite group that acts on $L$ 
preserving the bilinear form $(x.y)$. If $L^G \ne 0$, then the restriction of $(x.y)$ makes $L^G$
a lattice of discriminant dividing $(d |G|)^r$, where $r=\rk(L^G)$.
\ele
{\em Proof.} The $G$-module $L_\Q$ is semisimple, hence is a direct sum of
$G$-modules $L^G_\Q\oplus V$, where $V$ is a vector space over $\Q$ such that $V^G=0$.
If $x\in L^G_\Q$ and $y\in V$, then $(x.y)=(x.gy)$ for any $g\in G$. Since 
$\sum_{g\in G} gy\in V^G=0$, we have $(x.y)=0$. Thus $L_\Q=L^G_\Q\oplus V$
is an orthogonal direct sum. It follows that the discriminant of the restriction of the bilinear form
on $L$ to $L^G$ is non-zero, so $L^G$ is indeed a lattice. It is clear that the
finite abelian group $(L^G)^*/L^G$ is generated by at most $r$ elements. 
Thus it is enough to show that $(L^G)^*/L^G$ is annihilated by $d |G|$.

The map $L^G\to (L^G)^*$ is the composition of the natural maps
$$L^G\hookrightarrow L\lra L^*\lra (L^G)^*.$$
Since $L^G$ is a primitive sublattice of $L$, the last map here is surjective.
Thus any $a\in (L^G)^*$ is in the image of $L^*$, hence $d a$ is in the image of $L$.
Since $|G| a=\sum_{g\in G} ga$, we see that $(d |G|)a$ is in the image of $L^G$. $\Box$

\medskip

\ble \label{16may}
Let $\ell$ be a prime. Let $M$ be a free $\Z_\ell$-module of finite positive rank
equipped with a symmetric $\Z_\ell$-valued
bilinear form $(x.y)$. Let $\Ga$ be a group that acts on $M$ preserving the form $(x.y)$. 
If $L\subset M$ is a $\Z_\ell[\Ga]$-submodule 
such that the restriction of $(x.y)$ to $L$ has discriminant $d\not=0$, then 
$d\cdot (M/L)^\Ga$ belongs to the image of the natural map $M^\Ga\to (M/L)^\Ga$.

For any positive integer $n$ the image of the natural map 
$(M/\ell^n)^\Ga\to ((M/L)/\ell^n)^\Ga$ contains $d\cdot ((M/L)/\ell^n)^\Ga$.
\ele
{\em Proof.} Let $L^\perp\subset M$ be the orthogonal complement to $L$ with respect to
$(x.y)$. We have $L\cap L^\perp=0$ because $d\not=0$. Hence the natural map $L\oplus L^\perp\to M$
is injective.

Let $x\in M$. For $y\in L$
the map $y\mapsto (x.y)$ is an element of $\Hom_{\Z_\ell}(L,\Z_\ell)$.
Thus we get a map of $\Z_\ell$-modules $M \to \Hom_{\Z_\ell}(L,\Z_\ell)$.
Since $d\not=0$, the restriction of this map to $L$
is injective and has cokernel annihilated by $d$. Hence
there is a $z\in L$ such that $d(x.y)=(z.y)$ for all $y\in L$. Thus 
$d x-z\in L^\perp$, proving that $d M\subset L\oplus L^\perp$. We
summarise this in the following commutative diagram:
$$\xymatrix{d M \ \ar@{^{(}->}[r] \ar[d]& \ L\oplus L^\perp \ \ar[d]\ar@{^{(}->}[r]&M \ar[d]\\
d (M/L) \ \ar@{^{(}->}[r]& \ L^\perp \ \ar@{^{(}->}[r]&M/L}
$$
The group $\Ga$ preserves $L$ and $(x.y)$, hence $\Ga$ also preserves $L^\perp$; thus
all the arrows in the diagram are maps of $\Ga$-modules.
Since the homomorphism $L \oplus L^\perp \to L^\perp$ has a section,
the first claim of the lemma follows.

For $n\geq 1$ we obtain a commutative diagram of $\Z_\ell[\Ga]$-modules
\begin{equation}
\xymatrix{d M/\ell^n \ \ar[r] \ar[d]& \ L/\ell^n\oplus L^\perp/\ell^n \ 
\ar[d]\ar[r]&M/\ell^n \ar[d]\\
d (M/L)/\ell^n \ \ar[r]& \ L^\perp/\ell^n \ \ar[r]&(M/L)/\ell^n}
\label{w3}
\end{equation}
If $\alpha\in((M/L)/\ell^n)^\Ga$, then $d\alpha$ comes from 
$(L^\perp/\ell^n)^\Ga$.
Similarly to the previous case, the map $L^\perp/\ell^n\to (M/L)/\ell^n$ factors through
$M/\ell^n\to (M/L)/\ell^n$. This proves the lemma. $\Box$

\medskip

Let $\ell$ be a prime and let $N$ be a free $\Z_\ell$-module of finite positive rank.
The free $\Z_\ell$-module $\End_{\Z_\ell}(N)$ has a symmetric 
$\Z_\ell$-valued bilinear form $\Tr(xy)$, where $\Tr$ is the usual matrix trace.

Let $\Lambda\subset\End_{\Z_\ell}(N)$ be a $\Z_\ell$-subalgebra.
We write $\End_\Lambda(N)$ for the centraliser of $\Lambda$ in $\End_{\Z_\ell}(N)$,
that is, the set of $x\in \End_{\Z_\ell}(N)$ such that $x\lambda=\lambda x$
for all $\lambda\in \Lambda$. 

\ble \label{13aug}
If the restriction of the bilinear form $\Tr(xy)$ to $\End_\Lambda(N)$ has discriminant $d\not=0$, then
there is an integer $r\geq 0$ such that for all $n\geq 1$ we have
$$\ell^r\cdot\End_\Lambda(N/\ell^n) \subset\End_\Lambda(N)/\ell^n\subset \End_\Lambda(N/\ell^n).$$
\ele
{\em Proof.} Write $M=\End_{\Z_\ell}(N)$, $L=\End_\Lambda(N)$,
and let $L^\perp$ be the orthogonal complement to $L$ in $M$. 
Since $d\not=0$ we have $L\cap L^\perp=0$. In particular, the only element of $L^\perp$
commuting with $\Lambda$ is 0.

It is clear that $L$ and $L^\perp$ are saturated, free $\Z_\ell$-submodules of $M$.
Thus for all $n\geq 1$ we have $L/\ell^n\subset M/\ell^n$ and 
$L^\perp/\ell^n\subset M/\ell^n$.
Let $\ell^a$ be the highest power of $\ell$ dividing $d$ in $\Z_\ell$. 
We are in the situation of Lemma \ref{16may},
so we have commutative diagram~(\ref{w3}). The first row of (\ref{w3}) implies
$$
\ell^a\cdot(M/\ell^n)\subset L/\ell^n+L^\perp/\ell^n \subset M/\ell^n.
$$
By retaining only the elements commuting with $\Lambda$ we obtain
\begin{equation}
\ell^a\cdot\End_\Lambda(N/\ell^n)\subset \End_\Lambda(N)/\ell^n+
(L^\perp/\ell^n)\cap \End_\Lambda(N/\ell^n) \subset \End_\Lambda(N/\ell^n).
\label{w1}
\end{equation}

We claim that there exists a positive integer $b$ such that for all $n> b$ we have
\begin{equation}(L^\perp/\ell^n)\cap \End_\Lambda(N/\ell^n)\subset \ell \cdot(L^\perp/\ell^n).\label{au13}
\end{equation}
Indeed, let $S_n$ be the subset of the left hand side consisting of the elements that are
not contained in $\ell \cdot(L^\perp/\ell^n)$.
Reduction mod $\ell^n$ maps $S_{n+1}$ to $S_n$ for each $n\geq 1$.
If all the finite sets $S_n$ are non-empty, then 
$\varprojlim S_n\not=\emptyset$. Any $x\in \varprojlim S_n$ is an element of $L^\perp
\setminus\ell L^\perp$, hence $x\not=0$. But $x\in \End_\Lambda(N)=L$,
contradicting $L\cap L^\perp=0$.

From (\ref{au13}), in view of a canonical isomorphism 
$\ell \cdot(L^\perp/\ell^n)\tilde\lra L^\perp/\ell^{n-1}$, for each $n>b$ we obtain an injection
\begin{equation}(L^\perp/\ell^n)\cap \End_\Lambda(N/\ell^n)
\hookrightarrow (L^\perp/\ell^{n-1})\cap \End_\Lambda(N/\ell^{n-1}).\label{au13bis}
\end{equation}
This implies
\begin{equation}
\ell^b\cdot ((L^\perp/\ell^n)\cap \End_\Lambda(N/\ell^n))=0, \quad\quad n\geq 1. \label{w2}
\end{equation} 
Combining (\ref{w1}) with (\ref{w2}) proves the lemma with $r=a+b$. $\Box$

\subsection{Algebras}

Let $B$ be a separable semisimple algebra over a field $k$. Then $B$ is the product 
of matrix algebras $B_i=\Mat_{r_i}(D_i)$, where $D_i$ is a division $k$-algebra, for $i=1,\ldots,m$. 
Let $K_i$ be the centre of $D_i$ and let $d_i^2=\dim_{K_i}(D_i)$.
Here $K_i$ is a finite separable field extension of $k$.
We call the {\em intrinsic trace} of $x\in B$ the trace $\Tr_B(x)$ of the linear transformation of $B$
defined by the left multiplication by $x$. Write
$x=x_1+\ldots+x_m$, where $x_i\in B_i$.
The {\em relative reduced trace} 
$\tr_{B/k}:B\to k$ is defined as the sum of compositions of the usual reduced trace 
$\tr_{B_i/K_i}:B_i\to K_i$ of the central simple
$K_i$-algebra $B_i$ with the trace of the finite separable field extension $\Tr_{K_i/k}:K_i\to k$,
for $i=1,\ldots,m$, see \cite[Def. 9.13]{Rei03}. Thus
$$\tr_{B/k}(x)= \sum_{i=1}^m \tr_{B_i/k}(x_i)=\sum_{i=1}^m \Tr_{K_i/k}\tr_{B_i/K_i}(x_i).$$
These two natural notions of trace are related as follows, see \cite{Rei03}, formula (9.22):
\begin{equation}\Tr_B(x)=\sum_{i=1}^m d_i r_i\, \tr_{B_i/k}(x_i).\label{a2}
\end{equation}
The two notions of trace give rise to two symmetric bilinear forms on $B$ with values in $k$:

(1) The form $\tr_{B/k}(xy)$. This form is non-degenerate, see \cite[Thm.~9.26]{Rei03}. 

(2) The {\em intrinsic} bilinear form $\Tr_B(xy)$. 

Now let $k=\Q$. Let $\Lambda$ be an order in the semisimple $\Q$-algebra $B$. 
In other words, $\Lambda$ is a subring of $B$ such that $\Lambda\otimes_\Z\Q=B$. 
The restriction of $\tr_{B/\Q}$ to $\Lambda$ takes values in $\Z$ (see \cite[Thm.~10.1]{Rei03}),
so the bilinear form $\tr_{B/\Q}(xy)$ is integral on $\Lambda$.
We define the discriminant $\discr(\Lambda)$ to be the discriminant of the lattice $\Lambda$, equipped with this bilinear form.

If $k=\Q_\ell$, we similarly define the discriminant of an order in a semisimple $\Q_\ell$-algebra (well-defined up to multiplication by a square in $\Z_\ell^\times$).

The following two statements are undoubtedly well known.
For example, the implication ``$\ell \nmid \discr(\Lambda) \Rightarrow \Lambda/\ell$ is semisimple'' of Corollary~\ref{a1} is essentially \cite[Lemma~2.3]{MW95} (except that in \cite{MW95}, the discriminant is defined using the intrinsic trace $\Tr_B$, while we use the reduced trace $\tr_{B/\Q}$).
Nevertheless 
we give a detailed proof as we could not find the full statement of this proposition in the literature.

\bpr \label{b1}
Let $\ell$ be a prime and let $\Lambda$ be an order in a semisimple $\Q_\ell$-algebra. 
Then the following conditions are equivalent.

{\rm (i)} $\ell$ does not divide $\discr(\Lambda)$.

{\rm (ii)} for some positive integers
$n_1,\ldots,n_r$ we have $\Lambda\cong\oplus_{i=1}^r\Mat_{n_i}(O_{k_i})$, where
$O_{k_i}$ is the ring of integers of an unramified finite field extension $k_i/\Q_\ell$ for $i=1,\ldots,r$.

{\rm (iii)} the $\F_\ell$-algebra $\Lambda/\ell$ is semisimple.
\epr
{\em Proof.} By assumption $\Lambda$ is an order in the semisimple $\Q_\ell$-algebra
$B=\Lambda\otimes\Q_\ell$.

Let us first assume that this order is maximal.
Any maximal order $M\subset B$ is a direct sum of maximal orders of the simple
components of $B$, see \cite[Thm.~10.5 (i)]{Rei03}. This direct sum is an orthogonal direct sum 
for the bilinear form $\tr_{B/\Q_\ell}(xy)$, so it is enough to consider
a maximal order $M\subset\Mat_r(D)$, where $D$ is a division $\Q_\ell$-algebra.
By \cite[Thm.~12.8]{Rei03} there is a unique maximal order $O\subset D$; it is the integral closure of 
$\Z_\ell$ in $D$. By \cite[Thm.~17.3]{Rei03} any maximal
order in $\Mat_r(D)$ is conjugate to $\Mat_r(O)$ by an element of $\GL_r(D)$, 
so we have an isomorphism of rings $M\cong \Mat_r(O)$. From this
we get an $\F_\ell$-algebra isomorphism $M/\ell\cong \Mat_r(O/\ell)$.

It is well known that ${\rm rad}(O/\ell)=0$ if and only if 
$D$ is an unramified {\em field} extension of $\Q_\ell$ (see \cite[Thm.~14.3]{Rei03}). Then $O/\ell$ is a 
field extension of $\F_\ell$,
hence $\Mat_r(O/\ell)$ is a semisimple $\F_\ell$-algebra. This shows the equivalence of (ii) and (iii).

If $K$ is the centre of $D$, and $R$ is the integral closure of $\Z_\ell$ in $K$, then
we have (cf. Exercise 1 on p.~223 of \cite{Rei03})
\begin{equation}
\discr(\Mat_r(O))=\N_{K/\Q_\ell}(\discr(\Mat_r(O)/R))\cdot\discr(R/\Z_\ell)^{r^2\dim_K(D)}.
\label{form}
\end{equation}
This element of $\Z_\ell$ is not divisible by $\ell$ if and only if $D=K$ is an unramified
field extension of $\Q_\ell$, see \cite[Cor.~25.10]{Rei03}. This proves the equivalence of (i) and (ii).

Now suppose that $\Lambda$ is not a maximal order. Let us show that in this case
each of (i), (ii), (iii) is false. By \cite[Cor.~10.4]{Rei03}
there is a maximal order $M\subset B$ that contains $\Lambda$.
Since the index $[M:\Lambda]$ equals $\ell^a$ for an integer $a\geq 1$ and 
$\discr(\Lambda)=[M:\Lambda]^2\,\discr(M)$, we see that $\ell$ divides $\discr(\Lambda)$,
so (i) does not hold.

To show that (iii) does not hold we need to show that
${\rm rad}(\Lambda/\ell)\not=0$, for which it is enough to exhibit
a non-zero two-sided nilpotent ideal in $\Lambda/\ell$. Let $N=\Lambda\cap\ell M$.
This is a two-sided ideal in $\Lambda$, hence $N/\ell\Lambda$ is a two-sided ideal in
$\Lambda/\ell$. By Nakayama's lemma we have $M/(\Lambda+\ell M)\not=0$.
Since $\dim_{\F_\ell}(M/\ell)=\dim_{\F_\ell}(\Lambda/\ell)$,
the cardinalities of the kernel and the cokernel of the natural homomorphism 
$\Lambda/\ell\to M/\ell$ are equal, so
$N/\ell\Lambda\not=0$. This ideal of $\Lambda/\ell$ is nilpotent. Indeed,
take any $x\in N/\ell\Lambda$ and lift it to $\tilde x\in N\subset \ell M$.
Then 
$\tilde x^{a+1}\in\ell^{a+1}M\subset \ell\Lambda$,
hence $x^{a+1}=0$. 

This also implies that (ii) does not hold. Indeed, otherwise $\Lambda/\ell$ would be 
a semisimple $\F_\ell$-algebra, which it is not. $\Box$

\bco \label{a1}
Let $\Lambda$ be an order in a semisimple $\Q$-algebra. 
A prime $\ell$ does not divide $\discr(\Lambda)$ if and only if
the $\F_\ell$-algebra $\Lambda/\ell$ is semisimple.
\eco
{\em Proof.} Apply Proposition \ref{b1} to the order
$\Lambda\otimes_\Z\Z_\ell$ in the semisimple $\Q_\ell$-algebra $\Lambda\otimes\Q_\ell$. $\Box$

\subsection{Abelian varieties} \label{av}

Let $k$ be a field with a separable closure $\bar k$ and Galois group $\Ga_k=\Gal(\bar k/k)$.
Let $A$ be an abelian variety over $k$ and let $\ell$ be a prime different from char$(k)$.
For each positive integer $n$ the Kummer sequence gives rise 
to an exact sequence of $\Ga$-modules
\begin{equation} \label{uno}
0\lra \NS(\ov A)/\ell^n\stackrel{c_1}\lra\H^2_\et(\ov A,\mu_{\ell^n})\lra\Br(\ov A)[\ell^n]\lra 0
\end{equation}
Let $A^\vee$ be the dual abelian variety, and let 
$e_{\ell^n,A}:A[\ell^n]\times A^\vee[\ell^n]\to \mu_{\ell^n}$
be the Weil pairing. We have canonical isomorphisms of $\Ga$-modules
\begin{equation}
\H^2_\et(\ov A,\mu_{\ell^n})\cong\wedge^2\H^1_\et(\ov A,\mu_{\ell^n})(-1)\cong
(\wedge^2A^\vee[\ell]^n)(-1)\cong\Hom(\wedge^2A[\ell^n],\mu_{\ell^n})\label{quatro}
\end{equation}
and an injective map of $\Ga$-modules, cf. \cite[Section 3.3]{SZ08}:
\begin{equation} \label{tre}
\H^2_\et(\ov A,\mu_{\ell^n})\cong\Hom(\wedge^2A[\ell^n],\mu_{\ell^n})
\hookrightarrow \Hom(A[\ell^n],A^\vee[\ell^n]).
\end{equation}
Here the image consists of those $u:A[\ell^n]\to A^\vee[\ell^n]$ such that $e_{\ell^n,A}(x,ux)=0$
for all $x\in A[\ell^n]$, that is, the form $e_{\ell^n,A}(x,uy)$ is alternating.

Let $\Hom(A[\ell^n],A^\vee[\ell^n])_{\rm sym}$ be the subgroup of {\em symmetric}
(or self-dual) homomorphisms $u:A[\ell^n]\to A^\vee[\ell^n]$. It is shown in
\cite[Remark 3.2]{SZ08} that $u$ in $\Hom(A[\ell^n],A^\vee[\ell^n])$ is symmetric if and only if 
$e_{\ell^n,A}(x,uy)=-e_{\ell^n,A}(y,ux)$, that is, the form $e_{\ell^n,A}(x,uy)$ is
skew-symmetric. All alternating forms are skew-symmetric, so we get an injective map
\begin{equation}
\H^2_\et(\ov A,\mu_{\ell^n})\hookrightarrow \Hom(A[\ell^n],A^\vee[\ell^n])_{\rm sym}. \label{m1}
\end{equation}
For $\ell\not=2$ all skew-symmetric forms are alternating, so that (\ref{m1})
is an isomorphism.

It is well known that $\NS(\ov A)$ is canonically isomorphic to the 
subgroup of self-dual elements
$\Hom(\ov A,\ov A^\vee)_{\rm sym}\subset \Hom(\ov A,\ov A^\vee)$. 
This allows one to rewrite the cycle map as a map of $\Ga$-modules
\begin{equation}
\Hom(\ov A,\ov A^\vee)_{\rm sym}/\ell^n\lra \H^2_\et(\ov A,\mu_{\ell^n}).\label{m2}
\end{equation}

\ble \label{lem1}
Let $k$ be a field of characteristic $0$.
The composition of maps $(\ref{m2})$ and $(\ref{m1})$ is the negative of the natural map $\Hom(\ov A,\ov A^\vee)_{\rm sym} \to \Hom(A[\ell^n],A^\vee[\ell^n])_{\rm sym}$ given by the action of 
endomorphisms of $\ov A$ on $\ell^n$-torsion points.
\ele
{\em Proof.}
The claim is that the following diagram commutes:
$$ \xymatrix{
   \NS(\ov A)/\ell^n                        \ar[r]^-{c_1} \ar[d]^-{\cong}
 & \H^2_\et(\ov A,\mu_{\ell^n})             \ar[r]^-{\cong}
 & \Hom(\wedge^2A[\ell^n],\mu_{\ell^n})     \ar[d]
\\ \Hom(\ov A,\ov A^\vee)_{\rm sym}/\ell^n  \ar[r]
 & \Hom(A[\ell^n],A^\vee[\ell^n])_{\rm sym} \ar[r]^-{[-1]}
 & \Hom(A[\ell^n],A^\vee[\ell^n])_{\rm sym}.
} $$
The vertical arrow on the left is induced by the map
$ \NS(\ov A)\to \Hom(\ov A,\ov A^\vee)$ that sends
 $\mathcal{L}$ to $\phi_{\mathcal{L}}$, 
where $\phi_{\mathcal{L}}$ is the morphism $\ov A \to \ov A^\vee$ defined in \cite[Ch.~6, Cor.~4]{Mum74}.
The vertical arrow on the right sends the Weil pairing $e^{\mathcal{L}}_{\ell^n}$, which is defined by
$$e^{\mathcal{L}}_{\ell^n}(x, y) = e_{\ell^n,A}(x, \phi_{\mathcal{L}}(y)),$$
to the restriction of $\phi_{\mathcal{L}}$ to $\ell^n$-torsion subgroups.
Thus it suffices to prove that going along the top of the diagram sends $\mathcal{L}$ to $-e^{\mathcal{L}}_{\ell^n}$.

For the proof we can assume that $k$ is finitely generated over $\Q$. 
Choose an embedding $\bar k\hookrightarrow\C$ and extend the ground field from $\bar k$ to $\C$.
Let $A(\C) = V/\Lambda$, where $V\cong \C^g$ is the tangent space to $A$ at $0$,
and $\Lambda$ is a lattice in $V$.

According to the Appell--Humbert theorem \cite[p.~20]{Mum74}, any line bundle $\mathcal{L}_{\C}$ on $A(\C)$ can be written in the form $\mathcal{L}(H,\alpha)$ for some Hermitian form $H$ on $V$ such that $E = \Im\, H$ takes integer values on $\Lambda \times \Lambda$ (and some additional data $\alpha$ which are not relevant to us here).
The first Chern class of $\mathcal{L}$ is given by $E \in \Hom(\wedge^2 \Lambda, \Z) \cong \H^2(A(\C), \Z)$.
Thus the top line of the above diagram takes $\mathcal{L} \bmod \ell^n$ to $\exp(2\pi i \ell^n E) \in \Hom(\wedge^2 (\ell^{-n}\Lambda/\Lambda), \mu_{\ell^n})$.

As explained on \cite[Ch.~24, p.~237]{Mum74}, if 
$x,y \in A(\C)[\ell^n]$ lift to $\tilde x, \tilde y \in \ell^{-n}\Lambda$, then
$$ \exp(-2\pi i \ell^n E(\tilde x, \tilde y)) = e^{\mathcal{L}}_{\ell^n}(x, y). $$
This completes the proof. $\Box$

\medskip

For any abelian variety $A$, let $E_A$ and $H_A$ be the $\Ga$-modules which make the following sequences exact:
\begin{gather}
0\lra\End(\ov A\times\ov A^\vee)\otimes\Z_\ell\lra \End_{\Z_\ell}(T_\ell(A)\oplus T_\ell(A^\vee))
\lra E_A\lra 0, \label{EA}
\\
0\lra \Hom(\ov A,\ov A^\vee)\otimes\Z_\ell\lra \Hom(T_\ell(A),T_\ell(A^\vee))\lra H_A\lra 0. \label{HA}
\end{gather}
Note that the exact sequence (\ref{HA}) is a direct summand of (\ref{EA}).

Using (\ref{uno}), (\ref{tre}) and Lemma~\ref{lem1}, we have a commutative diagram of $\Gamma$-modules with exact rows:
\begin{equation} \label{diagram}
\vcenter{\xymatrix{
   0                                                     \ar[r]
 & \NS(\ov A)/\ell^n                                     \ar[r] \ar[d]
 & \H^2_\et(\ov A,\mu_{\ell^n})                          \ar[r] \ar[d]
 & \Br(\ov A)[\ell^n]                                    \ar[r] \ar[d]
 & 0
\\ 0                                                     \ar[r]
 & \Hom(\ov A,\ov A^\vee)/\ell^n                         \ar[r]^-{[-1]}
 & \Hom(A[\ell^n],A^\vee[\ell^n])                        \ar[r]
 & H_A/\ell^n                                            \ar[r]
 & 0
}}
\end{equation}

\ble \label{ker2}
Let $k$ be a field of characteristic $0$.
The kernel of the homomorphism $\Br(\ov A)[\ell^n] \to H_A/\ell^n$ at the right of \eqref{diagram} has exponent dividing $2$.
\ele
{\em Proof.}
Let $C_1$ and $C_2$ denote the cokernels of the left and central vertical arrows of~\eqref{diagram} respectively.
Since the central vertical arrow is injective, the snake lemma implies that $\ker(\Br(\ov A)[\ell^n] \to H_A/\ell^n)$ injects into $\ker(C_1 \to C_2)$.

If $\rho \in \Hom(\ov A, \ov A^\vee)/\ell^n$ maps to $0$ in $C_2$,
then the image of $\rho$ in $\Hom(A[\ell^n], A^\vee[\ell^n])$ lies in the image of $\H^2_\et(\ov A,\mu_{\ell^n})$ and hence is in $\Hom(A[\ell^n], A^\vee[\ell^n])_\sym$ by (\ref{m1}).
Choose any $\tilde\rho \in \Hom(\ov A, \ov A^\vee)$ lifting $\rho$.
Then $\tilde\rho + \tilde\rho^\vee$ is a symmetric lift of $2\rho$.
Since $\NS(\ov A)/\ell^n \to \Hom(\ov A,\ov A^\vee)_{\rm sym}/\ell^n$ is an isomorphism,
we conclude that $2\rho$ maps to $0$ in $C_1$.
This shows that $2 \cdot \ker(C_1 \to C_2) = 0$, which proves the lemma.
$\Box$

\section{Equivalence of variants of conjectures of Coleman and Shafarevich} \label{2}

Let $A$ be an abelian variety of dimension $g\geq 1$ over a field $k$.
Then $\End(A)$ is a free abelian group of positive rank at most equal to $4g^2$; as a ring, it
is an order in the finite-dimensional semisimple algebra $\End(A)_\Q$, 
see \cite[Ch.~19, Corollaries 1 and 3]{Mum74}. 
In Section \ref{1} we defined two integral symmetric bilinear forms on any order in $\End(A)_\Q$,
in particular, on $\End(A)$. The action of $\End(A)$ on $A$ by endomorphisms gives rise to a third
bilinear form. To fix notation we review all these forms here.
\begin{itemize}
\item The bilinear form $\tr(xy)$ on $\End(A)$, 
where $\tr:\End(A)_\Q\to\Q$ is the {\em reduced trace}. 
We call the discriminant of this form $\discr(A)$.
\item The {\em intrinsic} integral symmetric bilinear form on $\End(A)$ is $\Tr_{\End(A)}(xy)$,
where $\Tr_{\End(A)}(x)$ is the trace of the linear map $\End(A)\to\End(A)$
sending $z$ to $xz$.
We call the discriminant of this form $\Delta_A$.
\item For any $a\in\End(A)$ and $n\in\Z$ the degree of the endomorphism $[n]-a$ of $A$
is a monic polynomial in $n$ with integer coefficients \cite[Ch.~19, Thm.~4]{Mum74}.
Let $\Tr_A(a)\in \Z$ be the negative of the coefficient of $n^{2g-1}$ in this polynomial.
For any prime $\ell$ not equal to char$(k)$,
we have that $\Tr_A(a)$ is equal to the trace of the $\Z_\ell$-linear transformation of the 
$\ell$-adic Tate module of $A$ defined by $a$. 
We call the discriminant of this form $\delta_A$.
\end{itemize}

\ble \label{2.2}
Let $g$ be a positive integer. We have $\Delta_A\not=0$, $\delta_A\not=0$.
There exist positive real constants $c_g$ and $C_g$, depending only on $g$,
such that for any abelian variety $A$ of dimension $g$ over a field $k$
we have $$c_g \leq |\Delta_A|/|\discr(A)|\leq C_g, \quad c_g \leq |\delta_A|/|\discr(A)|\leq C_g.$$
\ele
{\em Proof.} If $A$ is a simple abelian variety, then 
$\End(A)_\Q$ is a division algebra over $\Q$. Let $K$ be the centre of $\End(A)_\Q$.
Write $e=[K:\Q]$ and $d^2=\dim_K \End(A)_\Q$.

Let $m$ be a positive integer.
By the proof of \cite[Ch.~19, Lemma]{Mum74}, any $\Ga_\Q$-invariant linear map
$\phi:\End(A^m)\otimes\ov \Q\to\ov \Q$ satisfying $\phi(xy)=\phi(yx)$ for all $x$ and~$y$
is a rational multiple of the reduced trace. This implies that
each of $\Tr_{\End(A^m)}$ and $\Tr_{A^m}$ is a non-zero rational multiple of the reduced trace 
$\tr$ on $\End(A^m)_\Q$. In particular, each form is non-degenerate.
By evaluating at the identity element of $\End(A^m)$ we obtain
$$\frac{\Tr_{\End(A^m)}(x)}{ed^2m^2}=\frac{\Tr_{A^m}(x)}{2gm}=
\frac{\tr(x)}{edm}.$$
Since $ed$ divides $2g$ by \cite[Ch.~19, Cor., p.~182]{Mum74}, we see that each of
$\Tr_{\End(A^m)}(x)$ and $\Tr_{A^m}(x)$ is an integral multiple of the reduced trace.

Now let $A_1,\ldots, A_n$ be simple, pairwise non-isogenous abelian varieties over $k$,
of dimension $\dim(A_i)=g_i$.
Let $B=\prod_{i=1}^n A_i^{m_i}$ for some positive integers $m_1,\ldots,m_n$.
Then $\End(B)$ is the product of rings $\End(A_i^{m_i})$, hence the matrix 
of each of the three forms on $\End(B)$ is the direct sum of $n$ diagonal blocks.
We deduce that 
$$\Delta_B=\discr(B)\cdot\prod_{i=1}^n (d_im_i)^{e_id_i^2m_i^2}, \quad
\delta_B=\discr(B)\cdot\prod_{i=1}^n \left(\frac{2g_i}{e_id_i}\right)^{e_id_i^2m_i^2}$$
This proves the lemma for $B$.

Finally, an arbitrary abelian variety $A$ over $k$ is isogenous to some $B=\prod_{i=1}^n A_i^{m_i}$,
where $A_1,\ldots, A_n$ are simple and pairwise non-isogenous abelian varieties over $k$. 
Then $\End(A)$ and $\End(B)$
are orders in $\End(A)_\Q\cong \End(B)_\Q$. In each of the three cases, 
the bilinear forms on $\End(A)_\Q$ and $\End(B)_\Q$
are compatible under this isomorphism. We have
$$[\End(A):\End(A)\cap \End(B)]^2\cdot\discr(A)=[\End(B):\End(A)\cap \End(B)]^2\cdot\discr(B).$$
The same formula holds for the discriminants of the two other forms.
Hence $\Delta_A/\Delta_B=\delta_A/\delta_B=\discr(A)/\discr(B)$, which proves the statement for $A$. $\Box$

\bpr \label{2.3}
The following statements are equivalent:

{\rm (i)} Coleman's conjecture about $\End(A)$;

{\rm (ii)} $\discr(A)$ is uniformly bounded for all abelian varieties $A$ of bounded dimension
defined over a number field of bounded degree;

{\rm (iii)} same as {\rm (ii)}, with $\discr(A)$ replaced by $\delta_A$;

{\rm (iv)} same as {\rm (ii)}, with $\discr(A)$ replaced by $\Delta_A$.
\epr
{\em Proof.} The equivalence of (ii), (iii) and (iv) was established in Lemma \ref{2.2}.
It is clear that (i) implies (iv). It remains to show that (ii) implies (i).

The ring $\End(A)$ is an order in the semisimple $\Q$-algebra $\End(A)_\Q$, which has
dimension at most $4g^2$. 
Since $\discr(A)$ is bounded, only finitely many semisimple $\Q$-algebras, up to isomorphism,
can be realised as $\End(A)_\Q$. Indeed, let $B$ be a semisimple $\Q$-algebra, with
simple components $B_i$ for $i=1,\ldots,n$, such that $\dim_\Q(B)$ is bounded
and $B$ contains an order of bounded discriminant. Then
each $B_i$ is a matrix algebra over a division algebra $D_i$ with centre $K_i$
such that $\dim_\Q(B_i)$ is bounded.
Using Proposition \ref{b1} and formula (\ref{form}),
we see that the discriminants of the fields $K_i$ are bounded, hence these fields 
belong to a fixed finite set of number fields. By the same proposition, the division $K_i$-algebra $D_i$
has bounded rank and ramification, so there are only finitely many isomorphism classes of such algebras.

By the structure theorem for maximal orders over Dedekind domains \cite[Thm. 21.6]{Rei03} 
and the Jordan--Zassenhaus theorem \cite[Thm. 26.4]{Rei03}, we know that there are only
finitely many maximal orders in $B$, up to conjugation by an element of $B^\times$, 
see \cite[Section 26, Exercise 8]{Rei03}.
Hence there are only finitely many isomorphism classes of maximal orders in $B$.
It follows that there are only finitely many isomorphism classes of
orders of bounded discriminant, so only finitely many rings can be realised as $\End(A)$. $\Box$

\bde \label{Silver}
Let $p$ be $0$ or a prime number. Define $d_p(g)=|\GL(2g,\F_3)|$ if $p\not=3$, and
$d_p(g)=|\GL(2g,\Z/4)|$ if $p=3$. Let us write $d(g)=d_0(g)$.
\ede

\bthe \label{22a1}
Coleman's conjecture about $\End(A)$ is equivalent to Coleman's conjecture 
about $\End(\ov A)$.
\ethe
{\em Proof.} Proposition \ref{2.3} can be applied over the ground field
$k$ as well as over $\bar k$. Thus to prove the theorem it is enough to show that
the uniform boundedness of $\delta_A$ is equivalent to the uniform boundedness of 
$\delta_{\ov A}$. 
It is clear from the definition of $\Tr_A$ that
for any $a\in\End(A)$ we have $\Tr_A(a)=\Tr_{\ov A}(a)$.
We note that $\End(A)=\End(\ov A)^\Ga$. 
By a result of Silverberg \cite[Thm.~2.4]{Sil92}, 
the cardinality of the image $G$ of $\Ga$ in the automorphism group of $\End(\ov A)$
is bounded by $d(g)$.
Thus assuming the boundedness of $\delta_{\ov A}$, 
the boundedness of $\delta_A$ follows from Lemma \ref{2.1}.

Conversely, let $A$ be an abelian variety of dimension $g$ defined over a number field of degree at most $e$.
By Silverberg's result, the boundedness of $\discr(A)$, where $A$ is considered over a number field
of degree at most $e\cdot d(g)$, implies the boundedness of $\discr(\ov A)$.
$\Box$

\bthe \label{Sh-eq}
The following conjectures are equivalent:

{\rm (i)} Shafarevich's conjecture about $\NS(\ov X)$;

{\rm (ii)} Shafarevich's conjecture about $\NS(\ov X)^\Ga$;

{\rm (iii)} Shafarevich's conjecture about $\Pic(X)$.
\ethe
{\em Proof.} Let $\discr(\NS(\ov X))$ be the discriminant of the bilinear form on $\NS(\ov X)$ given by
the intersection pairing. Define $\discr(\NS(\ov X)^\Ga)$ and $\discr(\Pic(X))$ similarly.
The ground field $k$ being of characteristic $0$, the ranks of these lattices do not exceed $20$.
By \cite[Ch.~9, Thm.~1.1]{cassels:quadratic-forms}, 
(i) is equivalent to the boundedness of $\discr(\NS(\ov X))$
for K3 surfaces defined over number fields of bounded degree, and similarly for (ii) and (iii).
It remains to prove the equivalence of these three boundedness conditions.

The boundedness of $\discr(\NS(\ov X))$ is equivalent to that of $\discr(\NS(\ov X)^\Ga)$
in view of Lemma \ref{2.1} and the classical Minkowski's lemma that gives a bound on the size of finite
subgroups of $\GL(n,\Z)$ in terms of $n$. To complete the proof it is enough to show that
$\Pic(X)$ is a subgroup of $\NS(\ov X)^\Ga$ of bounded index. 
The spectral sequence $\H^p(k,\H^q(\ov X,\G_m))\Rightarrow\H^{p+q}(X,\G_m)$ gives rise to 
the well known exact sequence of low degree terms
$$0\lra \Pic(X)\lra \Pic(\ov X)^\Ga\lra\Br(k)\to\Br(X).$$

The degree of the second Chern class of the tangent bundle of a K3 surface is $24$ \cite[Ch.~1, 2.4]{Huy16}, so every K3 surface has a $0$-cycle of degree $24$ (see the discussion of Chow groups in \cite[Ch.~12]{Huy16}).
This implies that there are
finite field extensions $k_1,\ldots,k_n$ of~$k$ such that $X$ has a $k_i$-point for each $i$, 
and ${\rm g.c.d.}([k_1:k],\ldots,[k_n:k])$ divides~$24$. If $K$ is a finite extension
of $k$ such that $X$ has a $K$-point, then the natural map $\Br(K)\to \Br(X_K)$ has a section
and so is injective. Now a restriction-corestriction argument shows that the kernel of
$\Br(k)\to\Br(X)$ is annihilated by $24$. It follows that $\Pic(X)$ is a subgroup of 
$\Pic(\ov X)^\Ga=\NS(\ov X)^\Ga$ such that the quotient group is annihilated by $24$.
Since $\NS(\ov X)^\Ga$ is a lattice of rank at most $20$, this implies that the index $[\NS(\ov X)^\Ga:\Pic(X)]$ is at most $24^{20}$. $\Box$

\section{Coleman implies Br(AV)} \label{3}

In order to prove that Coleman's conjecture implies Br(AV), we need to prove that Coleman implies two things for each pair of positive integers $(d,g)$:
\begin{enumerate}
\itemsep=0pt
\item[(1)] there exists a constant $C=C(d,g)$ such that for all primes $\ell > C$ and all abelian varieties of dimension $g$ defined over number fields of degree $d$, $\Br(\ov A)[\ell]^\Ga = 0$;
\item[(2)] for each prime $\ell \leq C$, there exists a constant $B(\ell,d,g)$ such that for all abelian varieties of dimension $g$ defined over number fields of degree $d$, $|\Br(\ov A)\{\ell\}^\Ga| \leq B(\ell,d,g)$.
\end{enumerate}

In~(2), $\Br(\ov A)\{\ell\}$ denotes the $\ell$-primary subgroup of $\Br(\ov A)$, that is, the set of elements whose order is a power of $\ell$.
Since $\Br(\ov A)^\Ga$ is a torsion abelian group, it is isomorphic to $\bigoplus_\ell \Br(\ov A)\{\ell\}^\Ga$ so (1) and (2) together suffice to prove Br(AV).

We prove the statement (1) at Theorem~\ref{3.2}(d) or Corollary~\ref{c2}, and (2) at Theorem~\ref{br-av-fixed-l}.
We apply Corollary~\ref{c2} and Theorem~\ref{br-av-fixed-l} to a set of abelian varieties containing a representative of every isomorphism class of abelian varieties of dimension~$g$ defined over a number field of degree~$d$.
Corollary~\ref{c2} and Theorem~\ref{br-av-fixed-l} show that, in order to prove Br(AV) for abelian varieties of dimension~$g$, it is sufficient to assume Coleman's conjecture (in the form of Proposition~\ref{2.3}(ii)) for abelian varieties of dimension at most $gQ(g)$ where we define
\[ Q(g)=[2g e^{2g/e}], \]
with $e$ denoting the base of the natural logarithm.
The significance of this quantity will appear in the proof of Theorem~\ref{autre1}.

\subsection{Abelian varieties at large primes, I} \label{Sect. 3.1}

We now show that Coleman's conjecture implies $\Br(\ov A)[\ell]^\Ga=0$
for abelian varieties $A$ of dimension $g$ defined over a number field of degree $d$, for all $\ell$
greater than some constant depending only on $d$ and $g$. 

\bthe \label{3.2}
Suppose that for all pairs of positive integers $(d,g)$ there is a constant $c=c(d,g)$
such that $|\discr(A)|<c$
for any abelian variety $A$ of dimension $g$
defined over a number field of degree $d$. Then there is a 
constant $C=C(d,g)$ such that for any prime $\ell>C$
and any abelian variety $A$ of dimension $g$
defined over a number field of degree $d$
we have the following statements.

{\rm (a)} the $\F_\ell$-algebra $\End(A)/\ell$ is semisimple;

{\rm (b)} the $\Ga$-module $A[\ell]$ is semisimple;

{\rm (c)} $\End(A)/\ell=\End(A[\ell])^\Ga$;


{\rm (d)} $\Br(\ov A)[\ell]^\Ga=0$.
\ethe

We give two proofs of Theorem~\ref{3.2}.
In this section we prove it via a shortcut provided by a recent theorem of R\'emond \cite[Thm.~1.1]{Rem18}.
In Section \ref{3.b} we prove a slightly stronger statement, which is valid over finitely generated fields of characteristic zero rather than just number fields, without using R\'emond's theorem.

Parts (a), (b) and~(c) of Theorem~\ref{3.2} can be proved by combining results of Masser and W\"ustholz \cite{MW95} with \cite[Thm.~1.1]{Rem18}.  The methods of Masser and W\"ustholz are similar to the proof given in this section.
See \cite[Lemma~2.3]{MW95} for part~(a), \cite[p.~222]{MW95} for part~(b) and \cite[Lemma~3.2]{MW95} for part~(c).

The result of R\'emond is used via the following lemma.

\ble \label{3.1}
Suppose that for all pairs of positive integers $(d,g)$ there is a constant $c=c(d,g)$
such that $|\discr(A)|<c$
for any abelian variety $A$ of dimension $g$
defined over a number field of degree $d$. Then for all pairs of positive integers $(d,g)$
there is a positive
integer $r=r(d,g)$ such that for any abelian variety $A$ of dimension $g$
defined over a number field of degree $d$, for any positive integer $n$
and any $\Ga$-submodule $W\subset A[n]$ there is an isogeny $u:A\to A$
such that $rW\subset uA[n]\subset W$.
\ele
{\em Proof.} Let $A$ be an abelian variety over a number field $k$
such that $[k:\Q]=d$ and $\dim(A)=g$. Under our assumptions,
by \cite[Thm.~1.1]{Rem18}
for any abelian variety $B$ defined over $k$ and $k$-isogenous to $A$
there is a $k$-isogeny $A\to B$ of degree bounded in terms of $d$ and $g$.
One deduces the existence of a positive integer $r=r(d,g)$ such that, for every pair of isogenous abelian varieties $A$ and $B$ of dimension $g$ defined over a number field $k$ of degree $d$, $[r]:A\to A$ factors through
some $k$-isogeny $A\to B$.
The rest of proof is identical to the proof of \cite[Cor.~5.4.1]{Zar85}. $\Box$

\medskip

\noindent {\em Proof of Theorem \ref{3.2}.} (a) For $\ell$ not dividing $\discr(A)$, the semisimplicity of 
the $\F_\ell$-algebra $\End(A)/\ell$ follows from Corollary \ref{a1}.

(b) We follow the proof of \cite[Cor.~5.4.3]{Zar85}.
Assume that $\ell$ does not divide $r(d,g)\discr(A)$, 
where $r(d,g)$ is as in Lemma \ref{3.1}.
To prove that $A[\ell]$ is a semisimple $\Ga$-module it is enough to show that for any 
$\Ga$-submodule $W\subset A[\ell]$ there is an idempotent $\pi\in \End(A)/\ell$
such that $W=\pi A[\ell]$. We apply Lemma \ref{3.1} with $n=\ell$ to obtain
an isogeny $u:A\to A$ such that $W=uA[\ell]$, where we used that $\ell$ and $r$ are coprime.
Since $\End(A)/\ell$ is semisimple by (a), we can write $u(\End(A)/\ell)=\pi(\End(A)/\ell)$
for some idempotent $\pi\in \End(A)/\ell$. Then $W=\pi A[\ell]$. 

(c) We follow the proof of \cite[Cor.~5.4.5]{Zar85} which refers to \cite[3.4]{Zar77}.
Assume that $\ell$ does not divide any of the integers $\discr(A)$, $r(d,g)$, $r(d,2g)$. 


Let $D$ be the centraliser of $\End(A)/\ell$ in $\End(A[\ell])\cong\Mat_{2g}(\F_\ell)$. Since 
$\End(A)/\ell$ is a semisimple $\F_\ell$-algebra by (a), the $\End(A)/\ell$-module $A[\ell]$ is semisimple,
hence $D$ is a semisimple $\F_\ell$-algebra. By the double centraliser theorem, the centraliser of $D$ in
$\End(A[\ell])$ is $\End(A)/\ell$. 

Take any $\varphi\in \End(A[\ell])^\Ga$. To prove that $\varphi\in \End(A)/\ell$ we need to show
that $\varphi$ commutes with $D$.
Applying Lemma \ref{3.1} to the graph of $\varphi$ in $A[\ell]^{\oplus 2}$ and using 
that $\ell$ does not divide $r(d,2g)$,
we write the graph of $\varphi$ as $u A[\ell]^{\oplus 2}$ for some 
$u\in \Mat_2(\End(A)/\ell)$. Let $p_i:A[\ell]^{\oplus 2}\to A[\ell]$ 
be the projector to the $i$-th summand, for $i=1,2$. 
Since $p_1u$ is surjective, for each $x \in A[\ell]$ we can write $x = p_1u(y)$ for some $y \in A[\ell]^{\oplus 2}$.
Then since $p_1u, p_2u : A[\ell]^{\oplus 2} \to A[\ell]$ are maps of $D$-modules and $\varphi p_1u = p_2u$, we have for all $d\in D$:
$$\varphi(dx) = \varphi p_1u(dy) = p_2u(dy) = d.p_2u(y) = d\varphi(x).$$
This proves (c).
%

(d) This follows from (b), (c) applied to $A \times A^\vee$ and the following lemma.

\ble \label{4.1e}
Let $k$ be a field of characteristic $0$.
Let $A$ be an abelian variety over $k$ of dimension $g \geq 1$.
If $\ell > 4g$, the $\Ga$-module $A[\ell]$ is semisimple, and
$$ \End(A \times A^\vee)/\ell = \End_\Ga(A[\ell] \oplus A^\vee[\ell]), $$
then $\Br(\ov A)[\ell]^\Ga = 0$.
\ele
{\em Proof.}
Since the $\Ga$-module $A[\ell]$ is semisimple, so is 
$A^\vee[\ell] \cong \Hom(A[\ell],\mu_\ell)$.
Hence by Serre's theorem \cite{Ser94} and the fact that $\ell > 4g$, we deduce that $\End(A[\ell] \oplus A^\vee[\ell])$ is a semisimple $\Ga$-module.

Let $E_A$ and $H_A$ be the $\Ga$-modules from the exact sequences (\ref{EA}) and (\ref{HA}).
The assumption
$ \End(A \times A^\vee)/\ell = \End_\Ga(A[\ell] \oplus A^\vee[\ell]), $
together with the semisimplicity of $\End(A[\ell] \oplus A^\vee[\ell])$, implies that $(E_A/\ell)^\Ga = 0$.
Since the sequence (\ref{HA}) is a direct summand of (\ref{EA}), we deduce that $(H_A/\ell)^\Ga = 0$.
Because $\ell$ is odd, combining this with Lemma~\ref{ker2} implies that $\Br(\ov A)[\ell]^\Ga = 0$.
$\Box$

\medskip

This finishes the proof of Theorem \ref{3.2}. $\Box$

\medskip

\noindent{\bf Remark} The statement of Lemma \ref{4.1e} remains true in positive characteristic.
For the proof one has to replace the reference to Lemma~\ref{ker2} by a counting 
argument similar to the one used in \cite[Lemma 3.5]{SZ08}.

\subsection{Abelian varieties at large primes, II} \label{3.b}

The aim of this section is to prove a somewhat stronger version of Theorem \ref{3.2}, see
Corollaries \ref{c1} and \ref{c2}.

Let $A$ be an abelian variety over a field $k$. For a prime $\ell\not={\rm char}(k)$ let
$T_\ell(A)$ be the $\ell$-adic Tate module of $A$, and let
$\rho_{\ell,A}:\Ga\to \Aut_{\Z_\ell}(T_\ell(A))$ be the attached $\ell$-adic Galois representation.
We denote by $\Lambda_\ell(A)$ the $\Z_\ell$-subalgebra of $\End_{\Z_\ell}(T_\ell(A))$
generated by $\rho_{\ell,A}(\Ga)$. 
Write $V_\ell(A)=T_\ell(A)\otimes_{\Z_\ell}\Q_\ell$ and define
$$D_\ell(A)=\Lambda_\ell(A)\otimes_{\Z_\ell}\Q_\ell \subset \End_{\Q_\ell}(V_\ell(A)).$$ 
Thus $\Lambda_\ell(A)$ is an order in the $\Q_\ell$-algebra $D_\ell(A)$. Define
$$E(A)=\End(A)\otimes\Q, \quad E_\ell(A)=\End(A)\otimes\Q_\ell\subset \End_{\Q_\ell}(V_\ell(A)).$$
It is well known that $E(A)$ is a semisimple $\Q$-algebra \cite{Mum74}, so that $E_\ell(A)$ is a semisimple
$\Q_\ell$-algebra. The $\Z$-algebra
$\End(A)$ is an order in $E(A)$, hence the $\Z_\ell$-algebra $\End(A)\otimes\Z_\ell$ is an order in $E_\ell(A)$.
It is clear that $\Lambda_\ell(A)$ and $\End(A)\otimes\Z_\ell$ are commuting subalgebras of 
$\End_{\Z_\ell}(T_\ell(A))$, and $D_\ell(A)$ and $E_\ell(A)$ are commuting subalgebras of 
$\End_{\Q_\ell}(V_\ell(A))$. By the work of Weil, Tate, Zarhin, Faltings, Mori on the Tate conjecture
\cite{Zar75, Zar76, Fal83, FaltingsM, Mor85}
it is known that if $k$ is finitely generated over its prime subfield, 
then $D_\ell(A)$ is a semisimple $\Q_\ell$-algebra and 
$$\End_\Ga(T_\ell(A))=\End_{\Lambda_\ell(A)}(T_\ell(A))=\End(A)\otimes\Z_\ell.$$
This implies
\begin{equation}
E_\ell(A)=\End_{D_\ell(A)}(V_\ell(A)), \quad D_\ell(A)=\End_{E_\ell(A)}(V_\ell(A)),
\label{autre3}
\end{equation} 
where the second identity follows from the first by the double centraliser theorem.

\bpr \label{autre6}
Let $k$ be a field finitely generated over $\Q$. Let $A$ be an abelian variety over $k$
of dimension $g\geq 1$. If $\ell$ is a prime not dividing $\discr(\Lambda_\ell(A))$, then
the $\Ga$-module $A[\ell]$ is semisimple, $\End(A)/\ell$ is a semisimple
$\F_\ell$-algebra, and $$\End_\Ga(A[\ell])=\End(A)/\ell.$$
\epr
{\em Proof.} Because $k$ is finitely generated over $\Q$,
$\Lambda_\ell(A)$ is an order in the {\em semisimple} $\Q_\ell$-algebra $D_\ell(A)$.
By Proposition \ref{b1} the $\F_\ell$-algebra $\Lambda_\ell(A)/\ell$ is semisimple, thus
$A[\ell]$ is a semisimple $\Lambda_\ell(A)/\ell$-module, hence also a semisimple $\Ga$-module.

Also by Proposition \ref{b1} we have an isomorphism 
$\Lambda_\ell(A)\cong\oplus_{i=1}^r\Mat_{n_i}(O_{k_i})$, where
$O_{k_i}$ is the ring of integers of an unramified field extension $k_i/\Q_\ell$ and $n_i$
is a positive integer, for $i=1,\ldots,r$. Write $\F_i=O_{k_i}/\ell$ for the residue field
of $O_{k_i}$.

Using the fact that $\Mat_{n_i}(O_{k_i})$ is Morita-equivalent to $O_{k_i}$, we obtain that
for each $i=1,\ldots, r$ there exists a free $O_{k_i}$-module $T_i$ of finite rank such that
$T_\ell(A)\cong \oplus_{i=1}^r T_i^{\oplus n_i}$, 
where the action of $\Lambda_\ell(A)$ on $T_i^{\oplus n_i}=T_i\otimes_{O_{k_i}}O_{k_i}^{\oplus n_i}$
is induced by the natural action of $\Mat_{n_i}(O_{k_i})$ on $O_{k_i}^{\oplus n_i}$.
Hence $\End(A)\otimes\Z_\ell$, being the centraliser of $\Lambda_\ell(A)$
in $\End_{\Z_\ell}(T_\ell(A))$, is equal to $\oplus_{i=1}^r\End_{O_{k_i}}(T_i)$.
Thus $\End(A)/\ell=\oplus_{i=1}^r\End_{\F_i}(T_i/\ell)$ is a semisimple $\F_\ell$-algebra. On the other hand,
the centraliser of $\Lambda_\ell(A)/\ell=\oplus_{i=1}^r\Mat_{n_i}(\F_i)$ in 
$\End_{\F_\ell}(A[\ell])=\End_{\F_\ell}(\oplus_{i=1}^r (T_i/\ell)^{\oplus n_i})$
is also equal to $\oplus_{i=1}^r\End_{\F_i}(T_i/\ell)$. This finishes the proof. $\Box$

\bpr \label{autre7}
Let $k$ be a field finitely generated over $\Q$. Let $A$ be an abelian variety over $k$
of dimension $g\geq 1$. If $\ell>4g$ is a prime that
does not divide $\discr(\Lambda_\ell(A\times A^\vee))$, then
$\Br(\ov A)[\ell]^\Ga=0$.
\epr
{\em Proof.} 
This follows from Proposition \ref{autre6} applied to $A \times A^\vee$ and Lemma \ref{4.1e}.
$\Box$

\medskip

The following theorem is the main result of this section.

\bthe \label{autre1}
Let $k$ be a field finitely generated over $\Q$. Let $A$ be an abelian variety over $k$
of dimension $g\geq 1$. There exists an abelian variety $B$ over $k$ which is $k$-isogenous
to an abelian subvariety of $A^{Q(g)}$ 
such that $\End(B)\otimes\Z_\ell$ is isomorphic to a matrix algebra over $\Lambda_\ell(A)^{\rm op}$.
In particular, $\discr(\Lambda_\ell(A))$ divides $\discr(B)$.
\ethe

\noindent{\em Proof.} Let us first prove the statement in the isotypic case,
i.e. when $A$ is a power of a simple abelian variety. Then $E(A)$ is a simple $\Q$-algebra. 

Let us fix an embedding $\bar k\hookrightarrow\C$. The natural action of $E(A)$ on $\H_1(A_\C,\Q)$
gives rise to an embedding $E(A)\subset \End_\Q(\H_1(A_\C,\Q))$, so that $E(A)$ is a simple
$\Q$-subalgebra of the matrix algebra $\End_\Q(\H_1(A_\C,\Q))$ containing its centre $\Q\,\Id$.
By \cite[Thm.~4.3.2]{Her68} the centraliser $D(A)=\End_{E(A)}(\H_1(A_\C,\Q))$ is also
a simple $\Q$-subalgebra of $\End_\Q(\H_1(A_\C,\Q))$. Moreover, by \cite[p.~105]{Her68} we have
$$\dim_\Q(D(A)) \cdot \dim_\Q(E(A)) =\dim_\Q(\End_\Q(\H_1(A_\C,\Q)))=4g^2.$$
Next, $D(A)$ is isomorphic to
a matrix algebra over a division $\Q$-algebra $F$, say $D(A)\cong\Mat_m(F)$. 
Comparing dimensions over $\Q$ we see that $m$ divides $2g$.

Let $M\cong F^{\oplus m}$
be a simple left $D(A)$-module (unique up to isomorphism). Any left $D(A)$-module that has finite dimension over $F$
is isomorphic to a direct sum of finitely many copies of $M$; in particular, the left $D(A)$-module $D(A)$
is isomorphic to $M^{\oplus m}$ and the $D(A)$-module $\H_1(A_\C,\Q)$
is isomorphic to $M^{\oplus r}$, where $r$ divides $2g$. We obtain an isomorphism of left $D(A)$-modules
\begin{equation}\H_1(A^m_\C,\Q)=\H_1(A_\C,\Q)^{\oplus m}\cong M^{\oplus mr}\cong
D(A)^{\oplus r}. \label{autre2}
\end{equation}

The Tate module $T_\ell(A)$ is isomorphic to $\H_1(A_\C,\Z)\otimes\Z_\ell$ as an $\End(A) \otimes \Z_\ell$-module.
Hence $V_\ell(A)\cong \H_1(A_\C,\Q)\otimes\Q_\ell$ as an $E_\ell(A)$-module.
But $T_\ell(A)$ is naturally a Galois module; in fact, we know that $V_\ell(A)$ is a $D_\ell(A)$-module satisfying (\ref{autre3}).
From this and the definition of $D(A)$ it follows that $D_\ell(A)=D(A)\otimes_\Q\Q_\ell$.
Now (\ref{autre2}) gives rise to isomorphisms of left $D_\ell(A)$-modules
(hence also of $\Ga$-modules)
\begin{equation}V_\ell(A^m)=V_\ell(A)^{\oplus m}\cong D_\ell(A)^{\oplus r}. \label{autre4}
\end{equation}
Recall that $\Lambda_\ell(A)$ is an order, hence a lattice in $D_\ell(A)$.
Let $S$ be the lattice in $V_\ell(A^m)$ obtained from the lattice 
$\Lambda_\ell(A)^{\oplus r}\subset D_\ell(A)^{\oplus r}$ via the isomorphism (\ref{autre4}).
It is clear that $S$ is stable under the action of the Galois group $\Ga$.
We note that $T_\ell(A^m)$ is also a $\Ga$-stable lattice in $V_\ell(A^m)$. 
Hence for some positive integer $N$ we have $\ell^N S\subset T_\ell(A^m)$.
There exists an abelian variety $B$ over $k$ such that $S\cong T_\ell(B)$ as $\Ga$-modules
together with a $k$-isogeny $\alpha: B\to A^m$
of degree $|T_\ell(A^m)/\ell^N S|$ such that $\alpha_*(T_\ell(B))=\ell^N S$. From the construction
of $S$ we have
$$\End(B)\otimes\Z_\ell=\End_\Ga(T_\ell(B))=\End_\Ga(\Lambda_\ell(A)^{\oplus r})=
\End_{\Lambda_\ell(A)}(\Lambda_\ell(A)^{\oplus r})=\Mat_r(\Lambda_\ell(A)^{\rm op}).$$
Since $m$ divides $2g$, this finishes the proof in the isotypic case.

In the general case $A$ is isogenous to $\prod_{i=1}^s A_i$, where each $A_i$ is a power of a simple
abelian variety and $\Hom(A_i,A_j)=0$ for $i\not=j$. Fixing such an isogeny we obtain an
isomorphism of $\Ga$-modules $V_\ell(A)\cong \oplus_{i=1}^s V_\ell(A_i)$
and an isomorphism
of $\Q_\ell$-algebras $D_\ell(A)\cong \oplus_{i=1}^s D_\ell(A_i)$. 
For each $i=1,\ldots,s$ we construct an isomorphism of $\Ga$-modules
$V_\ell(A_i^{m_i})\cong D_\ell(A_i)^{\oplus r_i}$ as in (\ref{autre4}), where both $m_i$ and $r_i$
divide $2\,\dim(A_i)$.
Write $r=\prod_{i=1}^s r_i$. Then we have isomorphisms of $\Ga$-modules
$V_\ell(A_i^{m_ir/r_i})\cong D_\ell(A_i)^{\oplus r}$, which add up to an isomorphism
of $\Ga$-modules 
\begin{equation}
V_\ell(A')\cong D_\ell(A)^{\oplus r}, \quad \text{where}\quad A'=\prod_{i=1}^s A_i^{m_ir/r_i}.
\label{autre5}
\end{equation} 
For any $x>0$ we have $\log(x)\leq x/e$, whence we obtain 
$r\leq e^{\sum_{i=1}^s r_i/e}\leq e^{2g/e}$. Thus $A'$ is an abelian 
subvariety of $A^{Q(g)}$. 

Let $S$ be the lattice in $V_\ell(A')$ obtained from the lattice $\Lambda_\ell(A)^{\oplus r}
\subset D_\ell(A)^{\oplus r}$ via isomorphism (\ref{autre5}). As above,
there is an abelian variety $B$ over $k$ such that $S\cong T_\ell(B)$ as $\Ga$-modules
together with a $k$-isogeny $B\to A'$, for which there is an isomorphism
$\End(B)\otimes\Z_\ell\cong\Mat_r(\Lambda_\ell(A)^{\rm op})$. This proves the theorem. $\Box$

\medskip

\noindent{\bf Remark} By the Poincar\'e reducibility theorem there are only finitely many abelian subvarieties 
of a given abelian variety considered up to $k$-isogeny.
(In fact, the same is true up to $k$-isomorphism, see \cite{LOZ}.) When $k$ is finitely generated over $\Q$
each isogeny class of abelian varieties over $k$ consists of finitely many $k$-isomorphism classes. 
(The case of number fields was treated in \cite[Prop.~3.1]{Zar85}. For the case of arbitrary finitely 
generated fields see \cite[Thm.~2 and its proof on pp.~214--215]{FaltingsM}.)
Thus the abelian variety $B$ in Theorem \ref{autre1} belongs to a finite set of $k$-isomorphism classes
determined by $A$.

\vspace{-1pt}
\bco \label{c1}
Consider a set $\mathcal S$ of abelian varieties such that each $A\in \mathcal S$
is defined over a field $k_A$ finitely generated over $\Q$. Suppose that there is a constant $c$
such that for every $A \in \mathcal S$ and every abelian variety $B$ over $k_A$ which is $k_A$-isogenous 
to an abelian subvariety of $A^{Q(\dim(A))}$, we have $\discr(B)<c$.
Then for every prime $\ell>c$ and every $A\in \mathcal S$,
the $\Ga$-module $A[\ell]$ is semisimple, $\End(A)/\ell$ is a semisimple
$\F_\ell$-algebra, and $\End_\Ga(A[\ell])=\End(A)/\ell.$
\eco
\vspace{-5pt}
{\em Proof.} Combine Theorem \ref{autre1} with Proposition \ref{autre6}. $\Box$

\vspace{-1pt}
\bco \label{c2}
Consider a set $\mathcal S$ of abelian varieties such that each $A\in \mathcal S$
is defined over a field $k_A$ finitely generated over $\Q$. Suppose that there is a constant $c$
such that for every $A \in \mathcal S$ and every abelian variety $B$ over $k_A$ which is $k_A$-isogenous 
to an abelian subvariety of $A^{Q(2\dim(A))}$, we have $\discr(B)<c$. Then
for every prime $\ell>\max(c,4\dim(A))$ and every $A\in \mathcal S$ we have $\Br(\ov A)[\ell]^\Ga=0$.
\eco
\vspace{-5pt}
{\em Proof.} Combine Theorem \ref{autre1} with Proposition \ref{autre7}. $\Box$

\subsection{Abelian varieties at a fixed prime} \label{Sect. 3.2}

The following proposition develops \cite[Remark 5.4.7]{Zar85}.
For abelian varieties over number fields, this proposition can also be proved by combining \cite[Lemma~4.1]{MW95} with \cite[Thm.~1.1]{Rem18}.


\bpr \label{4.6}
Let $\ell$ be a prime and let $g$ be a positive integer.
Consider a set~$\mathcal S$ of abelian varieties such that each $A\in \mathcal S$
has dimension at most $g$ and
is defined over a field $k_A$ finitely generated over $\Q$. Suppose that there is a constant $c$
such that for every $A \in \mathcal S$ and every abelian variety $B$ over $k_A$ which is $k_A$-isogenous 
to an abelian subvariety of $A^{Q(\dim(A))}$, we have $\discr(B)<c$.

Then there exists a positive integer $a=a(\mathcal S)$ such that for every abelian variety
$A \in \mathcal S$ and every $n \geq 1$,
the subgroup 
$[\ell^a]\cdot \End_\Ga(A[\ell^{n+a}])$ of
$\End_\Ga(A[\ell^n])$ is contained in the image of $\End(A)/\ell^n$.
\epr
{\em Proof.} 
We can apply Lemma \ref{13aug} to the Tate module $N=T_\ell(A)$, where
$\Lambda=\Lambda_\ell(A)$ is the $\Z_\ell$-subalgebra of $\End_{\Z_\ell}(T_\ell(A))$
generated by $\rho_{\ell,A}(\Ga)$.
Indeed, by Faltings \cite{FaltingsM} we have $\End_\Lambda(N)=\End_\Ga(T_\ell(A))
=\End(A)\otimes\Z_\ell$, whereas the restriction of $(x.y)=\Tr(xy)$ to $\End(A)\otimes\Z_\ell$
is non-degenerate by Lemma \ref{2.2}. 

Lemma~\ref{13aug} now gives a positive integer~$a$ such that
$\ell^a \cdot \End_\Gamma(A[\ell^{n+a}])$
is contained in the image of $\End(A)$ for every $n \geq 1$.
However $a$ may depend on the subalgebra $\Lambda$ of $\End_{\Z_\ell}(T_\ell(A))$ and on the structure of $T_\ell(A)$ as a $\Lambda$-module.

As recalled in Section \ref{3.b}, $\Lambda$ is an order
in an semisimple $\Q_\ell$-algebra of dimension at most $16g^2$. 
By Theorem \ref{autre1}, $\discr(\Lambda)$ is bounded by $c$.
By a similar argument to that used in the proof of Proposition~{2.3}, there are only finitely many isomorphism classes of $\Z_\ell$-orders of given discriminant
in semisimple $\Q_\ell$-algebras of given dimension.
Thus there are finitely many possibilities for $\Lambda$.

Furthermore, for each $\Z_\ell$-algebra $\Lambda$, there are only finitely many isomorphism classes of $\Lambda$-modules of given finite $\Z_\ell$-rank.
This implies that our constant $a$ can be chosen to depend only on $\mathcal S$.
$\Box$

\medskip

The main result of this section is the following theorem.

\bthe \label{br-av-fixed-l}
Let $\ell$ be a prime and let $g$ be a positive integer.
Consider a set $\mathcal S$ of abelian varieties such that each $A\in \mathcal S$
has dimension at most $g$ and
is defined over a field $k_A$ finitely generated over $\Q$. Suppose that there is a constant $c$
such that for every $A \in \mathcal S$ and every abelian variety $B$ over $k_A$ which is $k_A$-isogenous 
to an abelian subvariety of $A^{Q(2\dim(A))}$, we have $\discr(B)<c$.

Then there
exists a positive integer $r$ such that for every abelian variety $A \in \mathcal S$
the group $\Br(\ov A)\{\ell\}^\Ga$ has cardinality dividing $\ell^r$.
\ethe
{\em Proof.}
Recall the definitions of $E_A$ and $H_A$ from the exact sequences (\ref{EA}) and (\ref{HA}).

We show first that there is an integer $s$ depending only on $\mathcal S$ such that, for every $n \geq 1$, we have $[\ell^s] \cdot (H_A/\ell^n)^\Ga = 0$.

We equip $\End_{\Z_\ell}(T_\ell(A)\oplus T_\ell(A^\vee))$ with the unimodular symmetric bilinear
form $\Tr(xy)$, where $\Tr$ is the usual matrix trace.
By Lemma~\ref{2.2} and our hypothesis on $\mathcal S$, the restriction of
this form to $\End(\ov A\times\ov A^\vee)\otimes_\Z\Z_\ell$ has bounded, non-zero discriminant.
By Lemma \ref{16may} this gives a positive integer $b$ such that $[\ell^b]\cdot (E_A/\ell^n)^\Ga$
is contained in the image of 
$$(\End_{\Z_\ell}(T_\ell(A)\oplus T_\ell(A^\vee))/\ell^n)^\Ga=\End_\Ga(A[\ell^n]\oplus A^\vee[\ell^n]),$$
for every $n\geq 1$.

By Proposition~\ref{4.6}, this implies that $[\ell^{a+b}] \cdot(E_A/\ell^n)^\Ga$ is contained in the image of $\End(A \times A^\vee)/\ell^n$.
But by the exact sequence (\ref{EA}), the image of $\End(A \times A^\vee)/\ell^n$ in $E_A/\ell^n$ is $0$.
Thus $[\ell^{a+b}] \cdot (E_A/\ell^n)^\Ga = 0$.
Since the exact sequence (\ref{HA}) is a direct summand of (\ref{EA}), it follows that $[\ell^{a+b}] \cdot (H_A/\ell^n)^\Ga = 0$, as required.

Hence by Lemma~\ref{ker2}, $[2\ell^{a+b}] \cdot \Br(\ov A)[\ell^n]^\Ga = 0$, so the exponent of $\Br(\ov A)\{\ell\}^\Ga$ divides $2\ell^{a+b}$.
Since $\Br(\ov A)[\ell]$ is a quotient of a free $\Z_\ell$-module $\H^2(\ov A,\Z_\ell(1))/(\NS(\ov A)\otimes\Z_\ell)$ of rank
\[ \phantom{\Box} \;\;
\rk(\H^2(\ov A,\Z_\ell)) - \rk(\NS(\ov A)) \leq g(2g-1)-1, \]
this implies a bound for the cardinality of $\Br(\ov A)\{\ell\}^\Ga$. $\Box$

\subsection{Converse results} \label{converse}

Let $p$ be $0$ or a prime number. The function $d_p(g)$ was introduced in Definition \ref{Silver}.

For an abelian group $B$ and a prime $p$, define $B(p')$ to be the subgroup of $B_{\rm tors}$ 
consisting of the elements whose order is not divisible by $p$.
For $p=0$, we write $B(p') = B_{\rm tors}$.

The following statement will be used in Section \ref{6}.

\bpr \label{3.3}
Let $k$ be a field of characteristic $p$.
Let $A$ be an abelian variety over $k$ of dimension $g\geq 1$.
If $\Br(\ov A\times \ov A^\vee)(p')^\Ga$ is annihilated by a positive integer~$M$,
then for any positive integer $n$ not divisible by $p$, we have
$$d_p(g)M\cdot \End_\Ga(A[n])\subset\End(A)/n \subset \End_\Ga(A[n]).$$
In particular,
if $\ell$ is a prime not dividing $d_p(g)M$ and not equal to $p$, then $$\End_\Ga(A[\ell])=\End(A)/\ell.$$
\epr
{\em Proof.} In \cite{SZ14} the last two authors used the Kummer sequence 
and the K\"unneth formula to obtain an expression for the Brauer group of a product
of varieties, see \cite[formula (20), p.~761]{SZ14}.
Applied to the abelian variety
$A\times A^\vee$ it gives a canonical isomorphism of $\Ga$-modules
$$\Br(\ov A\times \ov A^\vee)[n]\cong \Br(\ov A)[n]\oplus\Br(\ov A^\vee)[n]\oplus
\End_{\Z/n}(A[n])/\big(\End(\ov A)/n\big).$$
Thus $\left(\End_{\Z/n}(A[n])/\big(\End(\ov A)/n\big)\right)^\Ga$ is a subgroup
of $\Br(\ov A\times \ov A^\vee)[n]^\Ga$, and so is annihilated by $M$.  
In view of the exact sequence of $\Ga$-modules
$$0\lra \End(\ov A)/n\lra \End_{\Z/n}(A[n]) \lra \End_{\Z/n}(A[n])/\big(\End(\ov A)/n\big)\lra 0$$
we conclude that $M\cdot \End_\Ga(A[n])\subset (\End(\ov A)/n)^\Ga$.

Let $G$ be the image of $\Ga$ in $\Aut(\End(\ov A))$ via its natural action on $\End(\ov A)$.
By a result of Silverberg \cite[Thm.~2.4]{Sil92}, $G$ is a finite group of order dividing $d_p(g)$.
The exact sequence of $\Ga$-modules
$$0\lra \End(\ov A)\stackrel{[n]}\lra \End(\ov A) \lra \End(\ov A)/n\lra 0$$
comes from the same sequence considered as an exact sequence of $G$-modules.
It gives rise to the exact sequence of abelian groups
$$0\lra \End(A)/n\lra (\End(\ov A)/n)^\Ga\lra\H^1(G,\End(\ov A)),$$
where we took into account that $\End(\ov A)^G=\End(\ov A)^\Ga=\End(A)$.
Since $\H^1(G,\End(\ov A))$ is annihilated by $d_p(g)$,
we obtain $d_p(g)\cdot (\End(\ov A)/n)^\Ga\subset \End(A)/n$, and thus
$d_p(g)M\cdot (\End_\Ga(A[n])\subset \End(A)/n$. $\Box$

\medskip

We point out the following partial converse to Theorem \ref{3.2}.

\bco \label{4.4}
Let $\ell$ be a prime and let $k$ be a field of characteristic $p\not=\ell$.
Let $A$ be an abelian variety over $k$ of dimension $g\geq 1$.
If $\ell$ does not divide $d_p(g)$,
the $\Ga$-module $A[\ell]$ is semisimple and $\Br(\ov A\times \ov A^\vee)[\ell]^\Ga=0$,
then the following hold:

{\rm (a)} the $\F_\ell$-algebra $\End(A)/\ell$ is semisimple;

{\rm (b)} $\ell$ does not divide $\discr(A)$.
\eco
{\em Proof.} (a) Since the $\Ga$-module $A[\ell]$ is semisimple, the $\F_\ell$-algebra $\End_{\Ga}(A[\ell])$
is semisimple. By the second part of Proposition \ref{3.3} 
it coincides with $\End(A)/\ell$.

(b) This follows from (a) and Corollary \ref{a1}. $\Box$

\section{Coleman implies Shafarevich} \label{5}

In this section, we show that Coleman's conjecture implies Shafarevich's conjecture.
We use the Kuga--Satake construction to relate Hodge structures associated with K3 surfaces to abelian varieties.
In order to obtain a result which is independent of the degree of polarisation of the K3 surface, we use a K3 surfaces version of Zarhin's trick from \cite{OS} and \cite{She}, which is described in terms of orthogonal Shimura varieties.

We first recall how one constructs an orthogonal Shimura variety from a lattice~$L$ with signature $(2,n)$, $n\geq 1$.
Let $\bSO(L)$ be the group scheme over $\Z$ whose
functor of points associates to a ring $R$ the group $\SO(L\otimes_\Z R)$.
Let $\S={\rm Res}_{\C/\R}(\G_m)$ denote the Deligne torus and let $\Omega_L$ be the set of $h\in\Hom(\S
,\bSO(L)_\R)$ such that the associated
$\Z$-Hodge structure on $L$ is of K3 type, that is, the following properties are satisfied:
\begin{enumerate}
\item \( \dim((L\otimes_\Z\C)_h^{(1,-1)})= \dim((L\otimes_\Z\C)_h^{(-1,1)})= 1 \) and 
\( \dim((L\otimes_\Z\C)_h^{(0,0)}) = n \);
\item for every non-zero \( v \in (L\otimes_\Z\C)_h^{(1,-1)} \) 
we have \( (v, v) = 0 \) and \( (v, \bar{v}) > 0 \);
\item \( ((L\otimes_\Z\C)_h^{(1,-1)},(L\otimes_\Z\C)_h^{(0,0)}) = 0 \).
\end{enumerate}
Sending $h$ to $(L\otimes_\Z\C)_h^{(1,-1)}$ identifies
$\Omega_L$ with $\{[x]\in\P(L\otimes_\Z\C)\ |\ (x^2)=0, \  (x,\bar x)>0\}$,
which is a homogeneous space of $\SO(L\otimes_\Z \R)$.

Let $\K\subset \bSO(L)(\A_{\Q,\rm f})$ be a compact open subgroup.
The canonical model of the
associated Shimura variety $\Sh_\K(L):=\Sh_\K(\bSO(L)_\Q,\Omega_L)$ is a quasi-projective
variety over $\Q$. By construction, the $\C$-points of $\Sh_\K(L)$
parameterise $\Z$-Hodge structures on $L$ satisfying properties 1, 2, 3 above.

Suppose that $\K$ is neat, as defined in \cite[0.6]{Pin89}.
(As explained in \cite[p.~5]{Pin89}, every compact open subgroup of $\bSO(L)(\A_{\Q,\rm f})$ contains a neat open subgroup of finite index.)
This implies that $\K \cap \bSO(L)(\Q)$ is torsion-free.
Then for each prime~$\ell$ 
there is a lisse $\Z_\ell$-sheaf $L_\ell$ on $\Sh_\K(L)$
defined by the inverse system of
finite \'etale covers $\Sh_{\K(\ell^m)}(L)\to\Sh_\K(L)$,
where $\K(\ell^m)$ is the largest subgroup of $\K$ that acts trivially on $L/\ell^m$.
Thus, to a $k$-point $x$ of $\Sh_\K(L)$ there corresponds
a representation $\Gal(\ov k/k) \to\SO(L\otimes_\Z\Z_\ell)$.
Putting together these representations for all $\ell$ gives a representation
\begin{equation} \label{sv-rep}
\phi_x: \Gal(\ov k/k) \to \SO(L\otimes_\Z\hat\Z).
\end{equation}
This construction was also described in \cite[3.1]{CM} or \cite[2.2]{UY13}.

From a lattice $L$ with signature $(2,n)$, $n\geq 1$, one can also construct
a spin Shimura variety, see \cite[Section 3]{MP16}. Let $C(L)$ be
the Clifford algebra of $L$, and let $C^+(L)\subset C(L)$ be the even Clifford algebra.
Let ${\bf GSpin}(L)$ be the group $\Z$-scheme
whose functor of points associates to a ring $R$ the group of
invertible elements $g$ of $C^+(L\otimes_\Z R)$ such that 
$g(L\otimes_\Z R)g^{-1}=L\otimes_\Z R$. If
$\tilde\K\subset{\bf GSpin}(L)(\A_{\Q,\rm f})$ is a compact open subgroup,
we write $\Sh^{\rm spin}_{\tilde\K}(L)$
for the Shimura variety $\Sh_{\tilde\K}({\bf GSpin}(L)_\Q,\Omega_L)$.
Let $\K$ be the image of $\tilde\K$ in $\bSO(L)(\A_{\Q,\rm f})$.
By \cite[4.4]{And96} this subgroup of $\bSO(L)(\A_{\Q,\rm f})$ is compact and open.
The natural group homomorphism ${\bf GSpin}(L)_\Q \to \bSO(L)_\Q$ induces a morphism $\Sh^{\rm spin}_{\tilde\K}(L)\to \Sh_{\K}(L)$.
This morphism is finite and surjective by \cite[Thm.~2.4]{Orr13}, and
defined over $\Q$ because the Shimura datum \( ({\bf GSpin}(L)_\Q,\Omega_L) \) has reflex field $\Q$ \cite[App.~1]{And96}.

Let $\tilde\K_N\subset {\bf GSpin}(L)(\hat\Z)$ be the set of elements congruent
to 1 modulo $N$ in $C^+(L\otimes_\Z\hat \Z)$. If $\tilde \K\subset\tilde \K_N$ for $N \geq 3$,
then $\tilde \K$ and $\K$ are neat and
the morphism $\Sh^{\rm spin}_{\tilde\K}(L)\to \Sh_{\K}(L)$ is \'etale.
Rizov shows in \cite[Section~5.5, (32)]{Riz10} that this morphism
restricts to an isomorphism on each geometric connected component.
Thus $\Sh^{\rm spin}_{\tilde\K}(L)\to \Sh_{\K}(L)$ has a section defined
over a number field $E$ which only depends on $L$ and $\tilde\K$.

There is a finite morphism of Shimura varieties from $\Sh^{\rm spin}_{\tilde\K}(L)$ to a moduli space of abelian varieties, defined over $\Q$.
In order to construct this,
we find a skew-symmetric form on $C(L)$ following \cite[Ch.~4, 2.2]{Huy16}.
Indeed, we choose orthogonal elements $f_1, f_2 \in L$ satisfying $(f_1^2), (f_2^2) > 0$.
Then we can define a skew-symmetric form on $C(L)$ by $\pm\Tr_{C(L)}(f_1f_2v^*w)$, where $\Tr_{C(L)}$ is the intrinsic trace.
The action of ${\bf GSpin}(L)$ on this form is multiplication by the spinor norm
(see \cite[Ch.~4, Prop.~2.5]{Huy16} for proofs of these facts, as well as the correct choice of sign).
The group ${\bf GSpin}(L)$ injects into
the group of symplectic similitudes ${\bf GSp}(C(L))$ of this form. 

If $\tilde \K\subset\tilde \K_N$,
then we have a morphism from $\Sh^{\rm spin}_{\tilde\K}(L)$ to the Shimura variety 
$\Sh_{\Ga_N}({\bf GSp}(C(L))_\Q, {\mathcal H}^\pm)$, where $\Ga_N$ is the subgroup of 
${\bf GSp}(C(L))(\hat\Z)$ consisting of the elements that are congruent to 1 modulo $N$.
The latter Shimura variety is identified with 
the moduli variety ${\mathcal A}_{g,\delta,N}$ parameterising abelian varieties of
dimension $g=2^{n+1}$, polarisation type $\delta$ (explicitly computable in terms of $L$ and $f_1, f_2$)
and level structure of level $N$.
If $N \geq 3$, then ${\mathcal A}_{g,\delta,N}$ is a fine moduli space,
so we can define the {\em Kuga--Satake abelian scheme} $f:A\to\Sh_{\K}(L)_E$ as
the pullback of the universal family of abelian varieties on ${\mathcal A}_{g,\delta,N}$
to $\Sh^{\rm spin}_{\tilde\K}(L)$, and then, after extending the ground field from $\Q$ to $E$, to
$\Sh_{\K}(L)_E$.
(As above, $E$ is a number field over which there exists a section of 
$\Sh^{\rm spin}_{\tilde\K}(L)_E \to \Sh_{\K}(L)_E$.
The Kuga--Satake scheme depends on the choice of such a section.)

The left multiplication of $L\subset C(L)$ on $C(L)$ gives a homomorphism
$L\hookrightarrow\End_\Z(C(L))$ whose cokernel is torsion-free.
Since $C(L) = R^1 f_{{\rm an},*}\Z$ as sheaves on $\Sh_\K(L)_\C$,
this gives rise to a morphism of variations
of $\Z$-Hodge structures $$L\to C(L)\to \End_\Z(R^1f_{{\rm an},*}\Z).$$
Via the comparison theorems we get a morphism
of $\Z_\ell$-sheaves $L_\ell\to \End_{\Z_\ell}(R^1f_*\Z_\ell)$.

\bpr \label{p1}
Let $L$ be a unimodular lattice of signature $(2,n)$, $n\geq 1$.
Let $\tilde \K\subset\tilde \K_3\subset {\bf GSpin}(L)(\hat\Z)$ be a compact open subgroup
and let $\K$ be the image of $\tilde\K$ in $\bSO(L)(\hat\Z)$.
For a $\C$-point $s$ of $\Sh_{\K}(L)$,
write $L_s$ for the $\Z$-Hodge structure on $L$ parameterised by $s$.
Define $T_s$ to be the smallest primitive sub-$\Z$-Hodge structure of $L_s$ whose 
complexification contains $L_s^{(1,-1)}$.

If Coleman's conjecture about $\End(\ov A)$
holds for abelian varieties of dimension $2^{n+1}$, then the discriminant
of the restriction of the bilinear form on $L$ to $T_s$ is bounded by a constant 
that depends only on $n$ and $d$, provided that $s$ is defined over a number field of degree~$d$.
\epr
{\em Proof.} Define $N_s=L\cap (L_s\otimes_\Z\C)^{(0,0)}$.
Then $T_s$ is the orthogonal complement to $N_s$ in $L$.
Since $L$ is unimodular, we have $|\discr(N_s)|=|\discr(T_s)|$, so it is 
enough to prove that $|\discr(N_s)|$ is bounded.

We equip the $\Z$-algebra $\End_\Z(C(L))={\rm Mat}_{2^{n+2}}(\Z)$ with 
the unimodular bilinear form $\Tr_{C(L)}(xy)$, where $\Tr_{C(L)}$ is the usual matrix trace (which,
by definition, is the same as the reduced trace).
This form is compatible with the Hodge structure, because the Hodge parameter $h_{\End}: \S \to \GL(\End_\Z(C(L)) \otimes \R)$ is given by $h_{\End}(z)(x) = h_{C(L)}(z) \, x \, h_{C(L)}(z)^{-1}$.
From the definition of the Clifford algebra we see that
the restriction of this form to $L\subset \End_\Z(C(L))$ is the original unimodular form
on $L$ multiplied by $2^{n+2}$. Let $L^\perp$ be the orthogonal complement
to $L$ in $\End_\Z(C(L))$. By the non-degeneracy of the form on $L$
we have $L\cap L^\perp=0$.
The index of
$L\oplus L^\perp$ in $\End_\Z(C(L))$ is equal to the discriminant of the restriction to $L$ of the bilinear form on $\End_\Z(C(L))$, and so depends only on $n$.

Suppose $s$ is defined over a number field $k$.
Without loss of generality we can assume that $k$ contains $E$. Thus we have an abelian variety $A_s$
defined over $k$ which is the fibre of $f:A\to\Sh_{\K}(L)$
at $s$. The natural injection $\End(A_{s,\C})\hookrightarrow \End_\Z(\H^1(A_{s,\C},\Z))$
gives an identification 
$$\End(A_{s,\C})=\End_\Z(\H^1(A_{s,\C},\Z))\cap \End_\C(\H^1(A_{s,\C},\C))^{(0,0)}.$$
In particular, $\End(A_{s,\C})$ is saturated in $\End_\Z(\H^1(A_{s,\C},\Z))$.

Since $L_s$ is a sub-$\Z$-Hodge structure of $\End_\Z(\H^1(A_{s,\C},\Z))$ and the Hodge structure is compatible with the bilinear form on $\End_\Z(\H^1(A_{s,\C},\Z)) \cong \End_\Z(C(L))$, we see
that $L^\perp_s$ is also a sub-$\Z$-Hodge structure.
Write $N'_s=L^\perp\cap (L^\perp_s\otimes_\Z\C)^{(0,0)}$. 
In the category of $\Q$-Hodge structures,
$\End_\Q(\H^1(A_{s,\C},\Q))=L_s\otimes_\Z\Q\oplus L^\perp_s\otimes_\Z\Q$, and so
$$\End(A_{s,\C})\otimes_\Z\Q=(N_s\otimes_\Z\Q)\oplus (N'_s\otimes_\Z\Q).$$
It follows that $N_s\oplus N'_s$ has finite index in $\End(A_{s,\C})$.

We have $\End(A_{s,\C})\cap L=N_s$ and $\End(A_{s,\C})\cap L^\perp=N'_s$.
If $x\in L$ and $y\in L^\perp$ are such that $(x,y)\in (L\oplus L^\perp)\cap \End(A_{s,\C})$,
then both $x$ and $y$ have type $(0,0)$, hence $x\in N_s$ and $y\in N'_s$.
This shows that $(L\oplus L^\perp)\cap \End(A_{s,\C})=N_s\oplus N'_s$.
Consequently the index of $N_s\oplus N'_s$ in $\End(A_{s,\C})$ divides
the index of $L\oplus L^\perp$ in $\End_\Z(C(L))$. Hence $|\discr(N_s)|$
divides the product of $|\discr(\End(A_{s,\C}))|$ and a constant depending only on $n$.
By Proposition \ref{2.3} and Theorem \ref{22a1}, Coleman's conjecture implies that
$|\discr(\End(A_{s,\C}))|$ is bounded. This finishes the proof. $\Box$

\medskip

The construction of the orthogonal Shimura variety associated to a lattice 
of signature $(2,n)$ is functorial
with respect to primitive embeddings of such lattices $\iota:L\hookrightarrow L'$.
Indeed, $\iota$ induces an injective group homomorphism of algebraic groups
$\bSO(L)_\Q\to \bSO(L')_\Q$ (since $L'_\Q=\iota(L_\Q) \oplus L_\Q^\perp$) and thus an injective homomorphism
$r: \bSO(L)(\A_{\Q,\rm f}) \to \bSO(L')(\A_{\Q,\rm f})$.
If $\K \subset \bSO(L)(\A_{\Q,\rm f})$ and $\K' \subset \bSO(L')(\A_{\Q,\rm f})$
are compact open subgroups such that $r(\K)\subset \K'$, then this gives rise to a
finite morphism of $\Q$-varieties 
$f: \Sh_\K(L) \lra \Sh_{\K'}(L')$.
When $\K'$ is neat, this morphism is
compatible with the associated variations of Hodge structures on $L$ and $L'$,
as well as with the associated $\ell$-adic sheaves. In particular, a $\C$-point $x$ of $\Sh_\K(L)$ gives rise to an isometric embedding of
associated $\Z$-Hodge structures $L_x\to L'_{f(x)}$.

We apply these considerations to orthogonal Shimura varieties related to
moduli spaces of polarised K3 surfaces, giving
a version of Zarhin's trick for K3 surfaces as proposed in \cite{OS} and \cite{She}.
For a positive integer $d$ let $\Lambda_{2d}$ be the lattice
$E_8(-1)^{\oplus 2} \oplus U^{\oplus 2} \oplus \langle -2d \rangle$.
There exist a positive integer $n$ and a {\em unimodular} lattice $\Lambda_\#$ of signature $(2,n)$
such that for each $d\geq 1$ there is a primitive embedding 
$\Lambda_{2d} \to \Lambda_\#$. In the version of \cite{OS} 
this lattice has been chosen as the even lattice $E_8(-1)^{\oplus 3} \oplus U^{\oplus 2}$
(so that $n=26$),
using results of Nikulin. Here we follow 
a simpler version based on Lagrange's four squares theorem as in \cite[Lemma 3.3.1]{She} and set
$\Lambda_\#= E_8(-1)^{\oplus 2} \oplus U^{\oplus 2} \oplus \langle -1 \rangle^{\oplus 5}$
(so that $n=23$).
For each $d$ we pick a primitive embedding $\iota_d:\Lambda_{2d} \to \Lambda_\#$, inducing
$r_d: \bSO(\Lambda_{2d})(\A_{\Q,\rm f}) \to \bSO(\Lambda_\#)(\A_{\Q,\rm f})$.
If $\K\subset \bSO(\Lambda_{2d})(\A_{\Q,\rm f})$ and 
$\K_\# \subset \bSO(\Lambda_\#)(\A_{\Q,\rm f})$ are
compact open subgroups such that
$r_d(\K)\subset \K_\#$, then there is a finite morphism of Shimura varieties over $\Q$
$$f_d: \Sh_{\K}(\Lambda_{2d}) \lra \Sh_{\K_\#}(\Lambda_\#).$$

Let $M_{2d}$ be the coarse moduli space over $\Q$ of primitively 
polarised K3 surfaces of degree $2d$;
this is a quasi-projective variety defined over $\Q$, 
see \cite[Ch.~5]{Huy16}. Let $\tilde{M}_{2d}$ be
the coarse moduli space over $\Q$ of triples $(X, \lambda, u)$ such that
$X$ is a K3 surface over a field of characteristic $0$,
$\lambda$ is a primitive polarisation of 
$X$ of degree~$2d$, and $u$ is an isometry
$$\det(P^2(\ov X, \Z_2(1))) \lra \det(\Lambda_{2d} \otimes_\Z \Z_2),$$
where $P^2(\ov X, \Z_2(1))$ is the orthogonal complement of the image of 
$\lambda$ in the $2$-adic \'etale cohomology $\H^2(\ov X, \Z_2(1))$.
We have a morphism $\tilde M_{2d}\to M_{2d}$, which is a double cover when $d>1$ and an isomorphism when $d=1$ \cite[Lemma~5.2 and Remark~5.8]{Tae}.
By the work of Rizov and Madapusi Pera (based on the Torelli theorem), 
there is an open immersion
$\tilde M_{2d}\hookrightarrow\Sh_{\K_d}(\Lambda_{2d})$ defined over $\Q$, where 
$$ \K_d = \{ g \in \SO(\Lambda_{2d}\otimes_\Z\hat \Z) : g \text{ acts trivially on } \Lambda_{2d}^*/\Lambda_{2d} \}. $$
For a proof that this immersion is defined over $\Q$, see \cite[Cor.~5.4]{MP15} (see also \cite[Thm.~3.9.1]{Riz10} and \cite{Tae}).

\bthe \label{5.2}
Coleman's conjecture about $\End(\ov A)$ implies Shafarevich's conjecture about $\NS(\ov X)$.
\ethe
{\em Proof.} Let $k$ be a number field and let $X$ be a K3 surface defined over $k$.
Let $d$ be a positive integer such that $X$ has a polarisation of degree $2d$ over $k$.
Then $X$ gives rise to a $k$-point on $M_{2d}$.
Replacing $k$ by a quadratic extension,
we can assume that this point lifts to a $k$-point $x$ on 
$\tilde M_{2d}\subset \Sh_{\K_d}(\Lambda_{2d})$.

Let $\K_\# \subset \bSO(\Lambda_\#)(\hat\Z)$ be the image of
$\tilde \K_3 \subset {\bf GSpin}(\Lambda_\#)(\hat\Z)$.
Define $$\K_d' = r_d^{-1}(\K_\#) \cap \K_d.$$
By \cite[Ch.~14, Prop.~2.6]{Huy16}, we have
\( r_d(\K_d) \subset \bSO(\Lambda_\#)(\hat\Z) \).
Hence $[\K_d:\K_d'] \leq [\bSO(\Lambda_\#)(\hat\Z):\K_\#]$, 
that is, the index $[\K_d:\K_d']$ is uniformly bounded.
Thus replacing $k$ by an extension of uniformly bounded degree we can assume that 
$x$ lifts to a $k$-point $x'$ on $\Sh_{\K_d'}(\Lambda_{2d})$.

We need to show that $|\discr(\NS(\ov X))|$ is universally bounded when $[k:\Q]$ is bounded. 
Choose an embedding $\bar k\hookrightarrow\C$. 
We have $\NS(\ov X)=\NS(X_\C)$. Let $T(X_\C)\subset\H^2(X_\C,\Z(1))$
be the transcendental lattice of $X_\C$ defined as the orthogonal complement
to $\NS(X_\C)$ in $\H^2(X_\C,\Z(1))$ with respect to the bilinear form given by the
cup-product. Since this form is unimodular, we have 
$|\discr(\NS(X_\C))|=|\discr(T(X_\C))|$, so it is enough to bound $|\discr(T(X_\C))|$.


Let $s=f_d(x') \in \Sh_{\K_\#}(\Lambda_\#)$. The proof of \cite[Lemma 4.3]{OS} shows that $\iota_d: \Lambda_{2d} \to \Lambda_\#$ induces
an isometry $T(X_\C)\tilde\lra T_s$.
Finally Proposition~\ref{p1} tells us that $|\discr(T_s)|$ is bounded by a constant that depends only on $[k:\Q]$. $\Box$

\section{Br(AV) implies V\'arilly-Alvarado} \label{6}

The main result of this section is that uniform boundedness of $\Br(\ov A)^\Gamma$, for abelian varieties $A$
of bounded dimension over number fields of bounded degree, implies V\'arilly-Alvarado's conjecture.

Before proving this main result, we relate two Galois representations attached to a polarised K3 surface $(X, \lambda)$ defined over a number field $k$. Choose
an isometry $u:\det(P^2(\ov X, \Z_2(1)))\to\det(\Lambda_{2d} \otimes_\Z \Z_2)$.
After replacing $k$ by a quadratic extension we can assume that $\Ga$ acts trivially on 
$\det(P^2(\ov X, \Z_2(1)))$. By \cite[Cor.~3.3]{Sai12}
the quadratic character through which $\Ga$ acts on $\det(\H^2(\ov X,\Q_\ell(1)))$ does not
depend on~$\ell$. Thus $\Ga$ acts trivially on $\det(P^2(\ov X, \Z_\ell(1)))$ for all primes $\ell$, hence
the representation $\rho_X: \Ga \to {\rm O}(P^2(\ov X, \hat\Z(1)))$ attached to $X$
takes values in $\SO(P^2(\ov X, \hat\Z(1)))$.

The triple $(X,\lambda,u)$ defines a $k$-point $x$ in 
$\tilde M_{2d}\subset \Sh_{\K_d}(\Lambda_{2d})$.
Choose a neat compact open subgroup $\K_d' \subset \K_d$
and let $x'$ be a lift of $x$ in $\Sh_{\K_d'}(\Lambda_{2d})$, defined over a finite extension $k'$ of $k$.
Let $\Ga' = \Gal(\ov k/k')$ and let $\phi_{x'}: \Ga' \to \SO(\Lambda_{2d} \otimes_\Z \hat\Z)$ denote the monodromy representation associated with the point~$x'$, as defined at (\ref{sv-rep}).

\vspace{-3pt}

\ble \label{6.1}
The adelic Galois representations $\rho_{X|\Ga'}: \Ga' \to \SO(P^2(\ov X, \hat\Z(1)))$ and $\phi_{x'}: \Ga' \to \SO(\Lambda_{2d} \otimes_\Z\hat\Z)$ are isometric.
\ele
{\em Proof.} This is an immediate consequence of \cite[Prop.~5.6(1)]{MP16}. $\Box$


\bthe \label{6.2}
Assume that for every positive integer $e$, there exists $B=B(e)>0$ such that every abelian variety $A$ of dimension $2^{25}$ defined over a number field of degree at most $e$ satisfies $|\Br(\ov A)^\Ga| < B$.
Then for every pair of positive integers $(e,M)$, there exists a constant $C=C(e,M)$ such that
for every K3 surface $X$ 
defined over a number field of degree~$e$, if $|\discr(\NS(\ov X))| < M$, 
then $|\Br(\ov X)^\Ga|<C$.
\ethe
{\em Proof.}
Let $X$ be a K3 surface defined over a number field $k$.
Let $d$ be a positive integer such that $X$ has a polarisation of degree $2d$ over $k$.
After an extension of the field $k$ of degree at most $2$,
$X$ is represented by a $k$-point $x$ of 
$\tilde M_{2d} \subset \Sh_{\K_d}(\Lambda_{2d})$.

Let $\Lambda_{2d}$, $\Lambda_\#$ and $\K_\#$, $\K_d'$ be the same as
in the proof of Theorem~\ref{5.2}. This proof shows that
replacing $k$ by an extension of uniformly bounded degree we can assume that 
$x$ lifts to a $k$-point $x'$ on $\Sh_{\K_d'}(\Lambda_{2d})$.
Let $s=f_d(x') \in \Sh_{K_\#}(\Lambda_\#)$. 

Write $\Lambda_{2d}\otimes_\Z\Z_\ell=\Lambda_{2d,\ell}$ and 
$\Lambda_{\#}\otimes_\Z\Z_\ell=\Lambda_{\#,\ell}$.
The injective homomorphism of $\Z$-modules $\iota_d:\Lambda_{2d}\to \Lambda_\#$ gives rise to
an injective homomorphism of $\Z_\ell$-modules $\iota_{d,\ell}:\Lambda_{2d,\ell}\to \Lambda_{\#,\ell}$,
which is also a homomorphism of $\Ga$-modules (with respect to the $\Ga$-module structures associated with the points $x' \in \Sh_{K_d'}(\Lambda_{2d})$ and $s \in \Sh_{K_\#}(\Lambda_\#)$ respectively).
Using comparison theorems between classical and \'etale cohomology, and noting that $T(X_\C) = \NS(X_\C)^\perp$, we see that
$T(X_\C)_\ell=T(X_\C)\otimes_\Z\Z_\ell$ has a canonical structure of a $\Ga$-module.
The proof of \cite[Lemma 4.3]{OS}, relying on \cite[Thm.~1.4.1]{Zar83}, shows that
$\iota_d(T(X_\C))=T_s$ (where $T_s$ is defined in Proposition~\ref{p1}), hence
$\iota_{d,\ell}$ sends $T(X_\C)_\ell$ isomorphically
onto $T_{s,\ell}=T_s\otimes_\Z\Z_\ell$. By Lemma~\ref{6.1}, we conclude that $\iota_{d,\ell}$ induces 
an isomorphism of $\Ga$-modules
$T(X_\C)_\ell\tilde\lra T_{s,\ell}.$

The Kummer exact sequence gives rise to short exact sequences of $\Gamma$-modules
\begin{equation} \label{kummer-seq}
0 \lra \NS(\ov X)/\ell^n \lra \H^2(\ov X,\mu_{\ell^n}) \lra \Br(\ov X)[\ell^n] \lra 0.
\end{equation}
Since the intersection pairing on $\H^2(X_\C,\Z(1))$ is unimodular, and using again the comparison between Betti and \'etale cohomology, this implies that there is an isomorphism of $\Ga$-modules
\begin{equation} \label{br-hom-t}
\Br(\ov X)[\ell^n] \cong \Hom(T(X_\C)_\ell, \Z/\ell^n).
\end{equation}

Since $\Br(\ov X)^\Ga$ is a torsion abelian group, in order to bound its size it suffices to show that $\Br(\ov X)[\ell]^\Ga = 0$ for large enough primes $\ell$ and to bound $|\Br(\ov X)\{\ell\}^\Ga|$ for every $\ell$.
Since $\Br(\ov X)\{\ell\}^\Ga \subset (\Q_\ell/\Z_\ell)^{20}$, in order to bound $|\Br(\ov X)\{\ell\}^\Ga|$ it suffices to bound the exponent of $\Br(\ov X)\{\ell\}^\Ga$.

Therefore using \eqref{br-hom-t}, in order to prove the theorem, it is enough to show that Br(AV) together with boundedness of $[k:\Q]$ and
$\discr(\NS(\ov X))$ imply that
\vspace{-2pt}
\begin{enumerate}
\itemsep=0pt
\item[(1)] there is a constant $C$ 
such that $\Hom_\Ga(T_{s,\ell},\Z/\ell)=0$ for any prime $\ell>C$,
where $s$ is any $k$-point of $\Sh_{\K}(\Lambda_{\#})$;

\item[(2)] for each prime $\ell$ there is an integer $m\geq 0$ 
such that $\ell^m\Hom_\Ga(T_{s,\ell},\Z/\ell^n)=0$ for any $n\geq 1$,
where $s$ is any $k$-point of $\Sh_{\K}(\Lambda_{\#})$.
\end{enumerate}

We assumed that $|\discr(\NS(\ov X))|=|\discr(T(X_\C))|=|\discr(T_s)| < M$,
thus the natural homomorphism of abelian groups $T_s\to \Hom_\Z(T_s,\Z)$
given by the intersection pairing is injective with cokernel of cardinality less than $M$.
Hence if $\ell \geq M$,
the $\Ga$-modules $T_{s,\ell}/\ell$ and $\Hom(T_{s,\ell},\Z/\ell)$
are canonically isomorphic, so to prove (1) it is enough to prove the following statement:
\vspace{-2pt}
\begin{enumerate}
\item[($1'$)] there is a constant $C$ such that $(T_{s,\ell}/\ell)^\Ga=0$ for any prime $\ell>C$.
\end{enumerate}
\vspace{-2pt}
For any fixed prime $\ell$ we have an injective homomorphism of $\Ga$-modules
$T_{s,\ell}\to\Hom(T_{s,\ell},\Z_\ell)$ with bounded cokernel. Thus, to prove (2)
it is enough to prove 
\vspace{-2pt}
\begin{enumerate}
\item[($2'$)] for each prime $\ell$ there is an integer $m\geq 0$ such that 
$[\ell^m]\cdot (T_{s,\ell}/\ell^n)^\Ga=0$ for all $n\geq 1$.
\end{enumerate}
\vspace{-2pt}

We use the notation of the proof of Proposition \ref{p1}.
Recall that $A=A_s$ is an abelian variety over $k$ of fixed dimension $g=2^{n+1}$
(where $\Lambda_\#$ has signature $(2,n)$ -- recall that we can take $n=23$).
We have an injective homomorphism of $\Z$-Hodge structures $T_s\to \Lambda_{\#,s} \to \End_{\Z}(\H_1(A_\C,\Z))$.
After tensoring with $\Z_\ell$ it gives rise to an
injective homomorphism of $\Ga$-modules $T_{s,\ell}\to \End_{\Z_\ell}(T_\ell(A))$.

We equip $\End_\Z(\H_1(A_\C,\Z))$ with 
the unimodular symmetric bilinear form $\Tr(xy)$, where $\Tr$ is the usual matrix trace.
After tensoring with $\Z_\ell$ this gives a $\Ga$-invariant form on $\End_{\Z_\ell}(T_\ell(A))$
with values in $\Z_\ell$.

Let $T_s^\perp$ be the orthogonal complement to $T_s$ in
$\End_{\Z}(\H_1(A_\C,\Z))$ with respect to $\Tr(xy)$. Clearly $T_s^\perp$ 
is saturated in $\End_{\Z}(\H_1(A_\C,\Z))$. In 
the proof of Proposition \ref{p1} we observed that the restriction of $\Tr(xy)$ to
$T_s$ is the intersection form on $T_s$ multiplied by $2^{n+2}$.
Since this form is non-degenerate, we have $T_s\cap T_s^\perp=0$.
The discriminant of $T_s$ is bounded by assumption and $\Tr(xy)$ is unimodular,
so
$$F=\End_{\Z}(\H_1(A_\C,\Z))/(T_s\oplus T_s^\perp)$$
is a finite abelian group of bounded size. We write $T_{s,\ell}^\perp=T_s^\perp\otimes_\Z\Z_\ell$.
This is the orthogonal complement to $T_{s,\ell}$ in $\End_{\Z_\ell}(T_\ell(A))$, so is naturally
a $\Ga$-module. In particular, $T_s/\ell^n$ and $T_s^\perp/\ell^n$ are $\Ga$-submodules of
$$\End_\Z(\H_1(A_{\C},\Z))/\ell^n=\End_{\Z_\ell}(T_\ell(A))/\ell^n=\End_{\F_\ell}(A[\ell^n])$$
for any prime $\ell$ and any positive integer $n$.

The bilinear form $\Tr(xy)$ is compatible with the Hodge structure.
Since $T_s\otimes\Q$ is an irreducible $\Q$-Hodge structure and contains elements of type $(1,-1)$,
it follows that all elements of $\End_{\Z}(\H_1(A_\C,\Z))$ of Hodge type $(0,0)$
are orthogonal to $T_s$.
In particular, $\End(A_\C)\subset T_s^\perp$.
Since $\End(A_\C)$ is saturated in $T_s^\perp \subset \End_\Z(H_1(A_\C,\Z))$, we also have $\End(A_\C)/\ell^n \subset T_s^\perp/\ell^n$.


We shall first prove ($2'$). Fix an arbitrary prime $\ell$ and let $\ell^a$ 
be the highest power of $\ell$ dividing the exponent of $F$.
Since $|F|$ is bounded, so is $a$.
Applying the snake lemma to the self-map $[\ell^n]$ of the exact sequence of $\Ga$-modules
$$0\lra T_{s,\ell}\oplus T_{s,\ell}^\perp\lra \End_{\Z_\ell}(T_\ell(A))\lra F[\ell^\infty]\lra 0$$
and then applying the left exact functor $-^\Gamma$, we get an exact sequence
\begin{equation} \label{ses}
0\lra F[\ell^n]^\Ga\lra (T_{s,\ell}/\ell^n)^\Ga\oplus (T_{s,\ell}^\perp/\ell^n)^\Ga\lra \End_\Ga(A[\ell^n]).
\end{equation}
By Proposition \ref{3.3} there is a positive integer $b$
that depends only on the upper bound for the cardinality of $\Br(\ov A\times \ov A^\vee)^\Ga$
such that $[\ell^b]\cdot \End_\Ga(A[\ell^n])$
is contained in $\End(A)/\ell^n \subset \End_\Ga(A[\ell^n])$. 
Recall that 
$$\End(A)/\ell^n=\End(\ov A)^\Ga/\ell^n\subset (\End(A_\C)/\ell^n)^\Ga\subset (T_{s,\ell}^\perp/\ell^n)^\Ga.$$
Let $x\in (T_{s,\ell}/\ell^n)^\Ga$. Using (\ref{ses}),
we see that there is a $y\in (T_{s,\ell}^\perp/\ell^n)^\Ga$ such that 
$\ell^b x\oplus y$ is in the image of $F[\ell^n]^\Ga$. 
Therefore, $\ell^a(\ell^bx\oplus y)=0$ 
hence $\ell^mx=0$, where $m=a+b$. This finishes the proof of ($2'$).

To prove ($1'$), note that if $\ell$ does not divide $|F|$, then $a=0$ in the above argument.
And by the second part of Proposition~\ref{3.3}, we can take $b=0$ for all primes $\ell$ greater than some constant depending only on the upper bound for
$\Br(\ov A\times \ov A^\vee)^\Ga$.
Thus for large enough primes $\ell$, the above argument for ($2'$) shows that $(T_{s,\ell}/\ell)^\Ga = 0$.
$\Box$

\vspace{-5pt} 
\appendix
\section{Coleman and Shafarevich's conjectures in one-parameter families}

In this appendix, we prove that Coleman and Shafarevich's conjectures for one-parameter families follow from a uniform open image theorem of Cadoret and Tamagawa.
Cadoret and Tamagawa's result concerns \( \ell \)-adic representations for a fixed \( \ell \);
we therefore also need a theorem of Hui on \( \ell \)-independence for the abelian varieties case, and the Mumford--Tate conjecture for the K3 surfaces case.

Throughout this appendix, we shall use the following notation.
Let \( k \) be a field finitely generated over \( \Q \).
Let \( X \) be a smooth geometrically connected variety over~\( k \) (in some of the theorems, \( X \) will be a curve).
Let \( \eta \) denote the geometric point and \( \ov\eta \) a geometric generic point of \( X \).
For a positive integer \( d \), we write
\[ \Xcld = \{ x \text{ a closed point of } X : [\kappa(x) : k] \leq d \}. \]
If \( \rho_\ell \colon \pi_1(X) \to \GL_m(\Z_\ell) \) is a representation, we write \( G_\ell = \rho_\ell(\pi_1(X)) \).
For any closed point \( x \) of \( X \), we write \( G_{x,\ell} = \rho_\ell \circ \sigma_x(\Gamma_{\kappa(x)}) \),
where \( \sigma_x \) denotes the map \( \Gamma_{\kappa(x)} \to \pi_1(X) \) induced by \( x \in X \).

\vspace{-3pt} 
\subsection*{Theorem of Cadoret and Tamagawa}
\vspace{-2pt} 

A representation \( \rho_\ell \colon \pi_1(X) \to \GL_m(\Z_\ell) \) is said to be \textit{GLP} if the Lie algebra of \( \rho_\ell(\pi_1(\ov X)) \) has trivial abelianisation.
Note that if a Lie algebra is semisimple, then its abelianisation is trivial, but the converse is not true.

We shall use the following example:
If \( Y \to X \) is a smooth proper scheme over \( X \), then the action of \( \pi_1(X) \) on the generic \( \ell \)-adic \'etale cohomology \( \H^{2i}(Y_{\ov\eta}, \Q_\ell(i)) \) is a GLP representation \cite[Thm.~5.8]{CT12}.
(See also Deligne's semisimplicity theorem for the monodromy of a smooth projective family of complex varieties \cite[Cor.~4.2.9]{Del71}).
This includes the case where \( A \to X \) is an abelian scheme and the representation is the action of \( \pi_1(X) \) on the generic Tate module \( T_\ell(A_\eta) \).


Our proof of Coleman and Shafarevich's conjectures for one-parameter families is based on the following theorem of Cadoret and Tamagawa.

\begin{theo} \textup{\cite[Thm~1.1]{CT13}} \label{ct}
Let \( k \) be a field finitely generated over \( \Q \).
Let \( X \) be a smooth geometrically connected curve over \( k \).
Let \( \rho_\ell \colon \pi_1(X) \to \GL_m(\Z_\ell) \) be a GLP representation.
Then for any integer \( d \geq 1 \), the set
\[ X_{\rho_\ell,d} = \{ x \in \Xcld : G_{x,\ell} \text{ is not open in } G_\ell \} \]
is finite.
\end{theo}

\subsection*{Coleman's conjecture}

We will prove the following theorem, of which Coleman's conjecture for one-parameter families is an immediate corollary.

\begin{theo} \label{av-finite}
Let \( k \) be a field finitely generated over \( \Q \).
Let \( X \) be a smooth geometrically connected curve over \( k \).
Let \( A \to X \) be an abelian scheme.
Then the set
\[ \{ x \in \Xcld : \End(\ov{A_x}) \neq \End(A_{\ov\eta}) \} \]
is finite.
\end{theo}

\begin{cor}
Let \( k \), \( X \), \( A \) be as in Theorem~\ref{av-finite}.
Then there are only finitely many isomorphism classes among the rings \( \End(\ov{A_x}) \), for \( x \in \Xcld \).
\end{cor}

The proof relies on the following result of Hui.

\begin{theo} \textup{\cite[Theorem~2.5]{Hui12}} \label{hui}
Let \( k \) be a field finitely generated over \( \Q \).
Let \( X \) be a smooth geometrically connected variety over \( k \).
Let \( A \to X \) be an abelian scheme.
For each prime \( \ell \),
let \( \rho_\ell \colon \pi_1(X) \to \Aut_{\Z_\ell}(T_\ell(A_\eta)) \) be the \( \ell \)-adic monodromy representation.
Then the set
\[ X_{\rho_\ell,d} = \{ x \in \Xcld : G_{\ell,x} \text{ is not open in } G_\ell \} \]
is independent of \( \ell \).
\end{theo}

We shall also need the following lemma.

\begin{lem} \label{av-ends}
Let \( k \) be a field finitely generated over \( \Q \).
Let \( X \) be a smooth geometrically connected variety over \( k \).
Let \( A \to X \) be an abelian scheme.
Let \( \rho_\ell \colon \pi_1(X) \to \Aut_{\Z_\ell}(T_\ell(A_\eta)) \) be the \( \ell \)-adic monodromy representation.
For any closed point \( x \) of \( X \), if \( G_{x,\ell} \) is open in \( G_\ell \), then
\[ \End(\ov{A_x}) \otimes \Z_\ell = \End(A_{\ov\eta}) \otimes \Z_\ell. \]
\end{lem}
\textit{Proof.}
Since the action of \( \Ga_{k(X)} \) on \( A_{\ov\eta}[3] \) is unramified, there is a finite \'etale cover \( X' \to X \) such that the \( 3 \)-torsion of \( A' = A \times_X X' \) is isomorphic (as a group scheme over \( X' \)) to the constant group scheme \( (\Z/3\Z)^{2g} \).

Choose a closed point \( x' \) of \( X' \) which maps to \( x \in X \).
Then all \( 3 \)-torsion of \( A'_{x'} = A_x \times_{\kappa(x)} \kappa(x') \) is defined over \( \kappa(x') \).
Hence by \cite[Thm.~2.4]{Sil92},
\[ \End(\ov{A_x}) = \End(A_x \times_{\kappa(x)} \kappa(x')). \]
Similarly,
\begin{equation} \label{EndAeta}
\End(A_{\ov\eta}) = \End(A_\eta \times_{k(X)} k(X')).
\end{equation}

Let
\[ G'_\ell = \rho_\ell(\pi_1(X')), \quad G'_{x,\ell} = \rho_\ell \circ \sigma_x(\Ga_{\kappa(x')}). \]
Now \( \kappa(x') \) is a finite extension of \( \kappa(x) \) and the restriction of \( \sigma_x \) to \( \Ga_{\kappa(x')} \) is the homomorphism \( \Ga_{\kappa(x')} \to \pi_1(X') \) induced by \( x' \in X' \).
Hence \( G'_{\ell,x} \subset G'_\ell \).
Furthermore \( G'_{x,\ell} \) is a closed, finite index subgroup of \( G_{x,\ell} \) so it is open in \( G_{x,\ell} \).
Hence the hypothesis that \( G_{x,\ell} \) is open in \( G_\ell \) implies that \( G'_{\ell,x} \) is open in \( G_\ell \).

Let $\alpha \colon \Ga_{k(X)} \to \pi_1(X)$ denote the natural
surjective homomorphism.  Let $L$ be the fixed field of
$(\rho_\ell \circ \alpha)^{-1}(G'_{x,\ell})$ inside $\ov{k(X)}$.
Since $G'_{x,\ell}$ is open in $G_\ell$, $L$ is a finite extension of
$k(X)$ and hence is a finitely generated field.
Furthermore \( \rho_\ell \circ \alpha(\Ga_L) = G'_{x,\ell} \).
Hence applying the Tate conjecture to the abelian variety \( A_\eta \times_{k(X)} L \), we get
\[ \End(A_\eta \times_{k(X)} L) \otimes \Z_\ell = \End_{\Ga_L}(T_\ell(A_\eta)) = \End_{G'_{x,\ell}}(T_\ell(A_\eta)). \]
Meanwhile the Tate conjecture for \( A_x \times_{\kappa(x)} \kappa(x') \) tells us that
\[ \End(A_x \times_{\kappa(x)} \kappa(x')) \otimes \Z_\ell = \End_{\Ga_{\kappa(x')}}(T_\ell(A_x)) = \End_{G'_{x,\ell}}(T_\ell(A_\eta)). \]
Because \( \alpha(\Ga_L) \subset \pi_1(X') \), we have \( k(X') \subset L \) and so \eqref{EndAeta} implies that
\[ \End(A_{\ov\eta}) = \End(A_\eta \times_{k(X)} L). \]
Combining the displayed equations proves the lemma.
\( \Box \)

\bigskip

\noindent \textit{Proof of Theorem~\ref{av-finite}.}
For each prime \( \ell \), let \( \rho_\ell \colon \pi_1(X) \to \Aut_{\Z_\ell}(T_\ell(A_\eta)) \) be the \( \ell \)-adic monodromy representation.
Let \( X_{\rho_\ell,d} \) be the set defined in Theorem~\ref{ct}.

By Theorem~\ref{ct}, \( X_{\rho_2,d} \) is finite.
So it suffices to show that \( \End(\ov{A_x}) = \End(A_{\ov\eta}) \) for all \( x \in \Xcld \setminus X_{\rho_2,d} \).

Consider a point \( x \in \Xcld \setminus X_{\rho_2,d} \).
By Theorem~\ref{hui}, \( x \not\in X_{\rho_\ell,d} \) for all~\( \ell \).
In other words, \( G_{x,\ell} \) is open in \( G_\ell \) for all~\( \ell \).
Hence by Lemma~\ref{av-ends},
\( \End(\ov{A_x}) \otimes \Z_\ell = \End(A_{\ov\eta}) \otimes \Z_\ell \)
for all~\( \ell \).
Since \( \End(A_{\ov\eta}) \) is a \( \Z \)-submodule of \( \End(\ov{A_x}) \), this implies that \( \End(A_{\ov\eta}) = \End(\ov{A_x}) \), as required.
\( \Box \)

\subsection*{Shafarevich's conjecture}

The proof of Shafarevich's conjecture for one-parameter families follows similar lines to that for Coleman's conjecture.
We will prove the following stronger result.

\begin{theo} \label{k3-finite}
Let \( k \) be a field finitely generated over \( \Q \).
Let \( X \) be a smooth geometrically connected curve over \( k \).
Let \( Y \to X \) be a smooth projective family of K3 surfaces.
Then the set
\[ \{ x \in X^{\cl,\leq d} : \NS(\ov{Y_x}) \neq \NS(Y_{\ov\eta}) \} \]
is finite.
\end{theo}

\begin{cor}
Let \( k \), \( X \), \( Y \) be as in Theorem~\ref{k3-finite}.
Then there are only finitely many isomorphism classes among the lattices \( \NS(\ov{Y_x}) \), for \( x \in \Xcld \).
\end{cor}

We prove the analogue of Hui's theorem using the Mumford--Tate conjecture for K3 surfaces, which was proved independently by Tankeev and Andr\'e \cite{Tan95,And96}.

\begin{prop}
Let \( k \) be a field finitely generated over \( \Q \).
Let \( X \) be a smooth geometrically connected variety over \( k \).
Let \( Y \to X \) be a smooth projective family of K3 surfaces.
For each prime \( \ell \),
let \( \rho_\ell \colon \pi_1(X) \to \Aut_{\Z_\ell}(\H^2(Y_{\ov\eta}, \Z_\ell(1))) \) be the \( \ell \)-adic monodromy representation.
Then the set
\[ X_{\rho_\ell,d} = \{ x \in \Xcld : G_{\ell,x} \text{ is not open in } G_\ell \} \]
is independent of \( \ell \).
\end{prop}
\textit{Proof.}
Since \( G_{\ell,x} \) and \( G_\ell \) are \( \ell \)-adic Lie groups, \( G_{\ell,x} \) is open in \( G_\ell \) if and only if \( \dim(G_{\ell,x}) = \dim(G_\ell) \).

Choose compatible embeddings \( k \hookrightarrow \C \) and \( \ov{k(X)} \hookrightarrow \C \).
By the Mumford--Tate conjecture, \( \dim(G_{\ell,x}) \) is equal to the dimension of the Mumford--Tate group of \( Y_x \times_k \C \) and hence is independent of \( \ell \).
Similarly \( \dim(G_\ell) = \dim(\MT(Y_{\ov\eta} \times_{\ov{k(X)}} \C)) \) so \( \dim(G_\ell) \) is independent of \( \ell \).
\( \Box \)

\medskip

We can now prove the following lemma and deduce Theorem~\ref{k3-finite} in the same way as Theorem~\ref{av-finite} is deduced from Lemma~\ref{av-ends}.

\begin{lem} \label{k3-ns}
Let \( k \) be a field finitely generated over \( \Q \).
Let \( X \) be a smooth geometrically connected variety over \( k \).
Let \( Y \to X \) be a family of K3 surfaces.
Let \( \rho_\ell \colon \pi_1(X) \to \Aut_{\Z_\ell}(\H^2(Y_{\ov\eta}, \Z_\ell(1))) \) be the \( \ell \)-adic monodromy representation.
For any closed point \( x \) of \( X \), if \( G_{x,\ell} \) is open in \( G_\ell \), then
\[ \NS(\ov{Y_x}) \otimes \Z_\ell = \NS(Y_{\ov\eta}) \otimes \Z_\ell. \]
\end{lem}
\textit{Proof.}
Let $\alpha \colon \Ga_{k(X)} \to \pi_1(X)$ denote the natural
surjective homomorphism.

Since the action of \( \Ga_{k(X)} \) on \( \NS(Y_{\ov\eta}) \) factors through a finite group, we can find a finite extension \( K' \) of \( k(X) \) such that this action becomes trivial after restricting to \( \Ga_{K'} \).
Similarly we can find a finite extension \( k_x' \) of \( \kappa(x) \) such that \( \Ga_{k_x'} \) acts trivially on \( \NS(\ov{Y_x}) \).
Since \( \alpha(\Ga_{K'}) \) is a closed, finite index subgroup of \( \pi_1(X) \), we may replace \( k'_x \) by a further finite extension such that \( \sigma_x(\Ga_{k'_x}) \subset \alpha(\Ga_{K'}) \) and \( \Ga_{k'_x} \) still acts trivially on \( \NS(\ov{Y_x}) \).

Let \( G'_{\ell,x} = \rho_\ell \circ \sigma_x(\Ga_{k_x'}) \).
This is a closed, finite index subgroup of \( G_{\ell,x} \), so it is open in \( G_{\ell,x} \).
Thus the hypothesis that \( G_{\ell,x} \) is open in \( G_\ell \) implies that \( G'_{\ell,x} \) is open in \( G_\ell \).

Let $L$ be the fixed field of
$(\rho_\ell \circ \alpha)^{-1}(G'_{x,\ell})$ inside $\ov{k(X)}$.
Since $G'_{x,\ell}$ is open in $G_\ell$, $L$ is a finite extension of
$k(X)$ and hence is a finitely generated field.
Furthermore \( \rho_\ell \circ \alpha(\Ga_L) = G'_{x,\ell} \).
Hence applying the Tate conjecture to the K3 surface \( Y_\eta \times_{k(X)} L \)
(see \cite{Tat94}, \cite{Tan95}, \cite{And96}), we get
\[ \NS(Y_{\ov\eta})^{\Ga_L} \otimes \Q_\ell = \H^2(Y_{\ov\eta}, \Q_\ell(1))^{\Ga_L} = \H^2(Y_{\ov\eta}, \Q_\ell(1))^{G'_{x,\ell}}. \]
Because \( \sigma_x(\Ga_{k'_x}) \subset \alpha(\Ga_{K'}) \), we have \( G'_{x,\ell} \subset \rho_\ell \circ \alpha(\Ga_{K'}) \) and so \( K' \subset L \).
Consequently \( \Ga_L \) acts trivially on \( \NS(Y_{\ov\eta}) \).
Since \( \NS(Y_{\ov\eta}) \otimes \Z_\ell \) is a saturated \( \Z_\ell \)-submodule of \( \H^2(Y_{\ov\eta}, \Z_\ell(1)) \) by the Kummer exact sequence~\eqref{kummer-seq}, we deduce that
\[ \NS(Y_{\ov\eta}) \otimes \Z_\ell = \H^2(Y_{\ov\eta}, \Z_\ell(1))^{G'_{x,\ell}}. \]

Meanwhile the Tate conjecture for \( Y_x \times_{\kappa(x)} k'_x \), together with a similar argument to the above to pass from \( \Q_\ell \)- to \( \Z_\ell \)-coefficients, tells us that
\[ \NS(\ov{Y_x}) \otimes \Z_\ell = \NS(\ov{Y_x})^{\Ga_{k'_x}} \otimes \Z_\ell = \H^2(\ov{Y_x}, \Z_\ell(1))^{\Ga_{k'_x}} = \H^2(Y_{\ov\eta}, \Z_\ell(1))^{G'_{x,\ell}}. \]
Combining the displayed equations proves the lemma.
\( \Box \)

\bigskip

\noindent Mathematics Institute, University of Warwick, Coventry CV4 7AL
England, U.K. 

\smallskip

\noindent {\tt martin.orr@warwick.ac.uk}

\bigskip

\noindent Department of Mathematics, South Kensington Campus,
Imperial College London, SW7 2BZ England, U.K. -- and --
Institute for the Information Transmission Problems,
Russian Academy of Sciences, 19 Bolshoi Karetnyi, Moscow 127994
Russia

\smallskip

\noindent {\tt a.skorobogatov@imperial.ac.uk}

\bigskip

\noindent Department of Mathematics, Pennsylvania State University, University Park, Pennsylvania 16802 USA 

\smallskip

\noindent {\tt zarhin@math.psu.edu}

\end{document}